\documentclass[11pt]{article}
\usepackage{amsfonts,amscd,amssymb}
\newtheorem{definition}{Definition}[section]
\newtheorem{lemma}[definition]{Lemma}
\newtheorem{proposition}[definition]{Proposition}
\newtheorem{corollary}[definition]{Corollary}
\newtheorem{remark}[definition]{Remark}

\newtheorem{theorem}[definition]{Theorem}

\def\va{\varepsilon}
\def\v{\varphi}

\def\rh{\rightharpoonup}
\def\lh{\leftharpoonup}

\def\ra{\rightarrow}
\def\a{\alpha}
\def\b{\beta}

\def\l{\lambda}
\def\r{\rho}
\def\cd{\cdot}

\def\O{\Omega}

\def\ov{\overline}

\def\mf{\mathfrak}
\def\mb{\mathbb}
\newcommand{\nat}{\mbox{$\;\natural \;$}}

\newcommand{\mfa}{\mbox{$\mf {a}$}}
\newcommand{\mfb}{\mbox{$\mf {b}$}}

\newcommand{\smi}{\mbox{$S^{-1}$}}
\newcommand{\gsm}{\mbox{$\blacktriangleright \hspace{-0.7mm}<$}}

\newcommand{\trl}{\mbox{${\;}$$\triangleright \hspace{-1.6mm}<$${\;}$}}
\newcommand{\gsl}{\mbox{${\;}$$>\hspace{-1.7mm}\triangleleft$${\;}$}}

\newcommand{\tx}{\mbox{$\tilde {x}$}}
\newcommand{\tX}{\mbox{$\tilde {X}$}}

\newcommand{\tqra}{\mbox{$\tilde {q}^1_{\rho }$}}
\newcommand{\tqrb}{\mbox{$\tilde {q}^2_{\rho }$}}

\renewcommand{\theequation}{\thesection.\arabic{equation}}

\def\rawo\lonra{\longrightarrow}

\def\ot{\otimes}

\newcommand{\eqref}[1]{(\ref{eq:#1})}

\newenvironment{proof}{{\it Proof.}}{\hfill $ \square $ \vskip 4mm}
\begin{document}
\title{On quasi-Hopf smash products and \\twisted tensor products 
of quasialgebras}
\author{Helena Albuquerque\thanks{The first author was partially supported 
by the Centre for Mathematics of the University of Coimbra (CMUC)} \\
Departamento de Matematica-Faculdade de Ciencias e Technologia,\\
Universidade de Coimbra, Apartado 3008, 3000 Coimbra, Portugal\\
e-mail: lena@mat.uc.pt 
\and Florin Panaite\thanks {Research
carried out while the second author was visiting the Centre for 
Mathematics of   
the University of Coimbra (CMUC) supported by a postdoctoral fellowship  
offered by CMUC. This author was also
partially supported by the programme CEEX of the Romanian 
Ministry of Education and
Research, contract nr. CEx05-D11-11/2005.}\\
Institute of Mathematics of the
Romanian Academy\\
PO-Box 1-764, RO-014700 Bucharest, Romania\\
 e-mail: Florin.Panaite@imar.ro}
\date{}
\maketitle
\begin{abstract}
We analyze some relations between quasi-Hopf smash products and 
certain twisted tensor products of quasialgebras. Along the way we 
obtain also some results of independent interest, such as a duality theorem 
for finite dimensional quasi-Hopf algebras  
and a universal property for generalized diagonal crossed products.
\end{abstract}
\section*{Introduction}
${\;\;\;}$If $(A, \mu _A, u_A)$ and $(B, \mu _B, u_B)$ are algebras in a 
monoidal category ${\cal C}$ and $R:B\otimes A\rightarrow A\otimes B$ 
is a morphism in ${\cal C}$ satisfying a certain list of axioms, then 
$A\otimes B$ becomes also an algebra in ${\cal C}$, with a 
multiplication defined in terms of $\mu _A$, $\mu _B$ and $R$. This 
construction appeared in a number of contexts and under different names. 
Following \cite{cap} we call such an $R$ a {\it twisting map} and the algebra  
structure on $A\otimes B$ afforded by it the {\it twisted tensor product} 
of $A$ and $B$ and denote it by $A\otimes _RB$ (if $A$ and $B$ are 
ordinary associative algebras and $R$ is the usual flip, then 
$A\otimes _RB$ coincides with the usual tensor product of algebras). 
Analogues of twisting maps for monads and operads are known  
as {\it distributive laws}, see for instance \cite{beck},  
\cite{markl}, \cite{street}. The twisted tensor product of associative 
algebras can be regarded as a representative for the cartesian product 
of noncommutative spaces, better suited than the ordinary tensor product, 
see \cite{cap}, \cite{jlpvo}, \cite{lpvo} for a detailed discussion and 
references.   
Prominent examples of twisted tensor products of algebras  
are the so-called {\it braided tensor products}, which are part of the 
''braided geometry'' developed by Majid in the early 1990's; namely, 
a braiding on a monoidal category provides a twisting map between any two 
algebras in the category. Hopf algebra theory provides also plenty of 
examples of twisted tensor products, in particular the usual smash product 
$A\# H$ is the twisted tensor product $A\otimes _RH$, with 
$R:H\otimes A\rightarrow A\otimes H$, $R(h\otimes a)=h_1\cdot a 
\otimes h_2$. \\
${\;\;\;}$The term {\it quasialgebra} was introduced in \cite{am1}, 
to designate an algebra in a monoidal category. We restrict here the term to 
algebras in monoidal categories associated to a quasi-bialgebra or quasi-Hopf  
algebra $H$, such as $_H{\cal M}$, $_H{\cal M}_H$, $_H^H{\cal YD}$. 
Since these categories have nontrivial associators, such quasialgebras are  
usually {\it nonassociative} as algebras, but their lack of 
associativity is (well) controlled by the associators of the category. 
Important examples of quasialgebras are the octonions and the other 
Cayley algebras, see \cite{am1}. Another class of examples of  
interest for us here is obtained as follows: if $H$ is a quasi-Hopf 
algebra, $B$ an associative algebra and $v:H\rightarrow B$ an algebra map, 
then on $B$ one can introduce a structure of algebra in $_H{\cal M}$, 
denoted by $B^v$, see \cite{bpv}. For the particular case $B=H$, $v=id_H$, 
the quasialgebra $B^v$ is denoted by $H_0$, and by \cite{bn} it 
follows that it is actually an algebra in the Yetter-Drinfeld category 
$_H^H{\cal YD}$. \\ 
${\;\;\;}$A feature of quasialgebras in this sense, not available for 
algebras in arbitrary monoidal categories (unless very special conditions 
hold, cf. \cite{schauenburg}) is the fact that there exist several  
{\it crossed products} associated to them, which are usually 
{\it associative} algebras. For instance, if $H$ is a quasi-bialgebra and 
$A$ is an algebra in $_H{\cal M}$, one can consider the smash product 
$A\# H$ introduced in \cite{bpv}, an associative algebra generalizing the 
classical Hopf smash product and sharing many of its properties. Other 
examples are diagonal crossed products, L-R-smash products,  
two-sided smash and crossed products, see \cite{bpvo}, \cite{hn1}, 
\cite{pv}.\\
${\;\;\;}$These quasialgebras are also part of an emerging 
{\it nonassociative geometry} (cf. \cite{akm}, \cite{bm}, 
\cite{majidgauge}), regarded as  
a further extension of noncommutative geometry, with the 
``coordinate algebra'' allowed to be nonassociative. By analogy with the 
associative situation, the twisted tensor product of quasialgebras 
might be regarded as a representative for the cartesian product of 
''nonassociative spaces''.\\
${\;\;\;}$The aim of this paper is to study certain classes of twisted tensor 
products of quasialgebras and their relations with quasi-Hopf smash products. 
If $H$ is a quasi-Hopf algebra, $A$ a left $H$-module algebra and $C$ an 
algebra in the Yetter-Drinfeld category $_H^H{\cal YD}$, we consider an 
object denoted by $A\diamond C$, which is a certain twisted tensor 
product in the category $_H{\cal M}$ and is defined in such a way that 
if $A$ is also an algebra in  $_H^H{\cal YD}$ then $A\diamond C$ coincides 
with the braided tensor product $A\underline{\otimes }C$ in $_H^H{\cal YD}$. 
We are mainly interested in the case $C=H_0$, so we have the left 
$H$-module algebra $A\diamond H_0$ and we want to see how is it  
related to the smash product $A\# H$. For Hopf algebras  
this question has a trivial answer: $A\# H$ and $A\diamond H_0$ coincide. 
This cannot be the case for proper quasi-Hopf algebras, because 
$A\# H$ is associative while $A\diamond H_0$ is {\it not}. The answer is the 
following: $A\diamond H_0$ is isomorphic as left $H$-module algebras 
to $(A\# H)^j$, where $j:H\rightarrow A\# H$, $j(h)=1\# h$. If $A$ 
is moreover an algebra in $_H^H{\cal YD}$, we prove that $(A\# H)^j$ 
becomes also an algebra in $_H^H{\cal YD}$ in a natural way and we have 
$A\underline{\otimes }H_0\simeq (A\# H)^j$ as algebras in $_H^H{\cal YD}$.\\
${\;\;\;}$More can be said if $A$ is $H_0$ itself. For this, we prove 
a result of independent interest: an algebra isomorphism $B^v\# H\simeq 
B\otimes H$, the quasi-Hopf analogue of a well-known result from Hopf  
algebra theory (see \cite{mon}).   
As a consequence, we obtain $H_0\underline{\otimes }H_0\simeq 
(H\otimes H)^{\Delta }$ as algebras in $_H^H{\cal YD}$. As another 
consequence, also of independent interest, we obtain a duality theorem, 
stating that if $H$ is a finite dimensional quasi-Hopf algebra then the 
two-sided crossed product $H\gsl H^*\trl H$ introduced by Hausser and Nill 
in \cite{hn1} is isomorphic to $End (H)\otimes H$ as algebras. \\
${\;\;\;}$In the last two sections we analyze some iterated products (in  
particular in the sense of \cite{jlpvo}) involving quasialgebras, and 
some universal properties. For instance, we find a relation between the 
universal property of $A\# H$ and the universal property of 
$A\diamond H_0$ regarded as a twisted tensor product of quasialgebras, 
and we find a new kind of universal property of the smash product $A\# H$, 
obtained as an immediate consequence of a very general universal 
property for generalized diagonal crossed products (as in \cite{bpvo}), 
which we obtain inspired by the universal property of the diagonal 
crossed product $H^*\bowtie {\mb A}$, formulated by Hausser and  
Nill in \cite{hn1}.                
\section{Preliminaries}\label{sec1}
\setcounter{equation}{0}
${\;\;\;}$
In this section we recall some definitions and results and we fix some 
notation that will be used throughout the paper.\\
${\;\;\;}$
We work over a base field $k$. All algebras, linear spaces 
etc. are over $k$; unadorned $\ot $ means $\ot_k$. Following 
Drinfeld \cite{d}, a quasi-bialgebra is a fourtuple $(H, \Delta ,
\va , \Phi )$, where $H$ is an associative algebra with unit $1$,  
$\Phi$ is an invertible element in $H\ot H\ot H$, and $\Delta :\
H\ra H\ot H$ and $\va :\ H\ra k$ are algebra homomorphisms
satisfying the identities (for all $h\in H$): 
\begin{eqnarray}
&&(id \ot \Delta )(\Delta (h))=%
\Phi (\Delta \ot id)(\Delta (h))\Phi ^{-1},\label{q1}\\
&&(id \ot \va )(\Delta (h))=h\ot 1, %
\mbox{${\;\;\;}$}%
(\va \ot id)(\Delta (h))=1\ot h,\label{q2} \\
&&(1\ot \Phi)(id\ot \Delta \ot id) (\Phi)(\Phi \ot 1)= (id\ot id
\ot \Delta )(\Phi ) (\Delta \ot id \ot id)(\Phi
),\label{q3}\\
&&(\va \ot id\ot id)(\Phi )=(id \ot \va \ot id )(\Phi )= 
(id \ot id\ot \va )(\Phi )=1\ot 1\ot 1.\label{q4}
\end{eqnarray}
The map $\Delta $ is called the coproduct or the
comultiplication, $\va $ the counit and $\Phi $ the reassociator.
We use the version of Sweedler's sigma notation: $\Delta (h)=h_1\ot  
h_2$, for $h\in H$, with the further convention 
$(\Delta \ot id)(\Delta (h))= h_{(1, 1)}\ot h_{(1, 2)}\ot h_2$ and     
$(id\ot \Delta )(\Delta (h))=h_1\ot h_{(2, 1)}\ot h_{(2,2)}$. 
We denote the tensor components of $\Phi$ 
by capital letters and those of $\Phi^{-1}$ by small letters: 
\begin{eqnarray*}
&&\Phi=X^1\ot X^2\ot X^3=T^1\ot T^2\ot T^3=Y^1\ot  
Y^2\ot Y^3=\cdots\\%
&&\Phi^{-1}=x^1\ot x^2\ot x^3=
t^1\ot t^2\ot t^3=y^1\ot y^2\ot y^3=\cdots
\end{eqnarray*}
${\;\;\;}$
The quasi-bialgebra $H$ is called a quasi-Hopf algebra if there exists an 
anti-automorphism $S$ of the algebra $H$ and elements $\a , \b \in
H$ such that, for all $h\in H$, we have:
\begin{eqnarray}
&&S(h_1)\a h_2=\va (h)\a \mbox{${\;\;\;}$ and ${\;\;\;}$}
h_1\b S(h_2)=\va (h)\b ,\label{q5}\\
&&X^1\b S(X^2)\a X^3=1 %
\mbox{${\;\;\;}$ and${\;\;\;}$}%
S(x^1)\a x^2\b S(x^3)=1.\label{q6}
\end{eqnarray}
The axioms for a quasi-Hopf algebra imply that $\va (\a )\va (\b 
)=1$, so, by rescaling $\a $ and $\b $, we may assume without loss
of generality that $\va (\a )=\va (\b )=1$ and $\va \circ S=\va $.\\
${\;\;\;}$
We recall that the definition of a quasi-bialgebra or  
quasi-Hopf algebra is 
"twist covariant" in the following sense. An invertible element
$F\in H\ot H$ is called a {\sl gauge transformation} or {\sl
twist} if $(\va \ot id)(F)=(id\ot \va)(F)=1$. If $H$ is a quasi-bialgebra 
or a quasi-Hopf algebra and $F=F^1\ot F^2\in H\ot H$ is a gauge  
transformation with inverse $F^{-1}=G^1\ot G^2$, then we can 
define a new quasi-bialgebra (respectively quasi-Hopf algebra) 
$H_F$ by keeping the 
multiplication, unit, counit (and antipode in the case of a quasi-Hopf 
algebra) of $H$ and replacing the reassociator, 
comultiplication and the elements $\alpha$ and $\beta$ by 
\begin{eqnarray}
&&\Phi_F=(1\ot F)(id \ot \Delta )(F) \Phi (\Delta \ot id)
(F^{-1})(F^{-1}\ot 1),\label{g2}\\
&&\Delta _F(h)=F\Delta (h)F^{-1},\;\;\;  \a_F=S(G^1)\a G^2,\;\;\;
\b_F=F^1\b S(F^2).\label{g3}
\end{eqnarray}
${\;\;\;}$
The antipode of a Hopf  
algebra is an anti-coalgebra morphism. For a quasi-Hopf algebra $H$, 
we have the following: there exists a gauge 
transformation $f\in H\ot H$ such that
\begin{equation} \label{ca}
f\Delta (S(h))f^{-1}=(S\ot S)(\Delta ^{cop}(h)) 
\mbox{,${\;\;\;}$for all $h\in H$.}
\end{equation}
The element $f$ can be computed explicitly. First set 
\begin{equation}
A^1\ot A^2\ot A^3\ot A^4=(\Phi \ot 1) (\Delta \ot id\ot id)(\Phi ^{-1}),
\end{equation}
\begin{equation} B^1\ot B^2\ot B^3\ot B^4=
(\Delta \ot id\ot id)(\Phi )(\Phi ^{-1}\ot 1),
\end{equation}
and then define $\gamma, \delta\in H\ot H$ by
\begin{equation} \label{gd}%
\gamma =S(A^2)\a A^3\ot S(A^1)\a A^4~~{\rm and}~~ \delta
=B^1\b S(B^4)\ot B^2\b S(B^3).
\end{equation}
Then $f$ and $f^{-1}$ are given by the formulae
\begin{eqnarray}
f&=&(S\ot S)(\Delta ^{cop}(x^1)) \gamma \Delta (x^2\b
S(x^3)),\label{f}\\%
f^{-1}&=&\Delta (S(x^1)\a x^2) \delta (S\ot S)(\Delta
^{cop}(x^3)).\label{g}
\end{eqnarray}
Moreover, $f$ satisfies the following relations:
\begin{eqnarray}
&&f\Delta (\alpha )=\gamma , \;\;\;\; \Delta (\beta )f^{-1}=\delta ,  
\label{moref}
\end{eqnarray}
and the corresponding twisted reassociator (cf. (\ref{g2})) is given by 
\begin{eqnarray}
&&\Phi _f=S(X^3)\otimes S(X^2)\otimes S(X^1). \label{muchmoref}
\end{eqnarray}
We record also the following relation from \cite{bn0}, for 
$f^{-1}=g^1\otimes g^2$ given by (\ref{g}):
\begin{eqnarray}
&&S^{-1}(\alpha g^2)g^1=\beta . \label{relg}
\end{eqnarray}
${\;\;\;}$
Suppose that $(H, \Delta , \varepsilon , \Phi )$ is a
quasi-bialgebra. If $U,V,W$ are left (right) $H$-modules, define
$a_{U,V,W}, {\bf a}_{U, V, W} :(U\otimes V)\otimes W\rightarrow
U\otimes (V\otimes W)$
by $a_{U,V,W}((u\otimes v)\otimes w)=\Phi \cdot (u\otimes
(v\otimes w))$ and 
${\bf a}_{U, V, W}((u\ot v)\ot w)= (u\ot (v\ot w))\cd \Phi ^{-1}$. 
The category $_H{\cal M}$ (${\cal M}_H$) of 
left (right) $H$-modules becomes a monoidal category (see
\cite{k}, \cite{m} for terminology) with tensor product 
$\otimes $ given via $\Delta $, associativity constraints
$a_{U,V,W}$ (${\bf a}_{U, V, W}$), unit $k$ as a trivial
$H$-module and the usual left and right
unit constraints.\\
${\;\;\;}$
Let again $H$ be a quasi-bialgebra. We say that a $k$-vector space 
$A$ is a left $H$-module algebra if it is an algebra in the
monoidal category $_H{\cal M}$, that is $A$ has a multiplication
and a usual unit $1_A$ satisfying the 
following conditions: 
\begin{eqnarray}
&&(aa')a''=(X^1\cd a)[(X^2\cd a')(X^3\cd
a'')],\label{ma1}\\
&&h\cd (a a')=(h_1\cd a)(h_2\cd a'), \;\;\;
h\cd 1_A=\va (h)1_A,
\end{eqnarray}
for all $a, a', a''\in A$ and $h\in H$, where $h\ot a\ra
h\cd a$ is the left $H$-module structure of $A$. Following
\cite{bpv} we define the smash product $A\# H$ as follows: as
vector space $A\# H$ is $A\ot H$ (elements $a\ot h$ will be
written $a\# h$) with multiplication
given by 
\begin{equation}\label{sm1}
(a\# h)(a'\# h')=%
(x^1\cd a)(x^2h_1\cd a')\# x^3h_2h'. 
\end{equation}
The smash product $A\# H$ is an 
associative algebra with unit $1_A\# 1_H$.\\ 
${\;\;\;}$
If $H$ is a quasi-Hopf algebra, $B$ an associative algebra and  
$v:H\rightarrow B$ an algebra map, then, following \cite{bpv}, we can 
introduce on the vector space $B$ a left $H$-module algebra structure, 
denoted by $B^v$ in what follows, for which the multiplication, unit and  
left $H$-action are:
\begin{eqnarray}
&&b\star b'=v(X^1)bv(S(x^1X^2)\alpha x^2X^3_1)b'v(S(x^3X^3_2)),\;\; 
\forall \;b, b'\in B, \label{mult}\\
&&1_{B^v}=v(\beta ), \;\;\;
h\triangleright _vb=v(h_1)bv(S(h_2)), \;\;\forall 
\;h\in H,\; b\in B.
\end{eqnarray}
${\;\;\;}$
If $H$ is a quasi-Hopf algebra and $A$ is a left $H$-module algebra, define 
the following maps: 
\begin{eqnarray}
&&j:H\rightarrow A\# H,\;\; j(h)=1\# h, \;\;\forall \; h\in H, \label{j}\\
&&i_0:A\rightarrow A\# H, \;\;i_0(a)=x^1\cdot a\# x^2\beta S(x^3), \;\;
\forall \; a\in A. \label{i0}
\end{eqnarray}
Then, by \cite{bpv}, $j$ is an algebra map and $i_0$ is a  
morphism of left $H$-module algebras from $A$ to $(A\# H)^j$. Moreover, 
the following universal property of the smash product $A\# H$ holds  
(see \cite{bpv}, Proposition 2.9): if $B$ is an associative algebra, 
$v:H\rightarrow B$ is an algebra map and $u:A\rightarrow B^v$ is a  
morphism of left $H$-module algebras, then there exists a unique algebra 
map $u\# v:A\# H\rightarrow B$ such that $(u\# v)\circ i_0=
u$ and $(u\# v)\circ j=v$; this map may be described 
explicitely as follows: 
\begin{eqnarray}
&&(u\# v)(a\# h)=v(X^1)u(a)v(S(X^2)\alpha X^3h), \;\;
\forall \; a\in A, \;h\in H.\label{unique}
\end{eqnarray}
We record the following relation from \cite{bpv}, which holds 
in $(A\# H)^j$ for all $a\in A$, $h\in H$:
\begin{eqnarray}
&&i_0(a)*j(h)=x^1\cdot a\# x^2hS(x^3). 
\label{inainte}
\end{eqnarray}
${\;\;\;}$
We recall now the invariance under twisting of the smash product (see for 
instance \cite{bn0}, \cite{bpv}, \cite{maj}).  
Let $H$ be a quasi-bialgebra, $F\in H\ot H$ a gauge transformation   
and $A$ a left $H$-module algebra. Then we can define a new  
multiplication on $A$, by 
$a\bullet a'=(G^1\cdot a)(G^2\cdot a')$, for all $a, a'\in A$,  
where $F^{-1}=G^1\otimes G^2$. 
If we denote by $A_{F^{-1}}$ the resulting structure, then 
$A_{F^{-1}}$ becomes a left $H_F$-module algebra, with the same unit and 
$H$-action as for $A$, and moreover the map  
$\pi: A\# H\rightarrow A_{F^{-1}}\# H_F,\;$   
$\pi (a\# h)=F^1\cdot a\# F^2h$, is an algebra isomorphism.\\ 
${\;\;\;}$
For further use we need also the notion of right $H$-module
algebra. Let $H$ be a quasi-bialgebra. We say that a $k$-linear
space $B$ is a right $H$-module algebra if $B$ is an algebra in
the monoidal category ${\cal M}_H$, i.e. $B$ has a multiplication
and a usual unit $1_B$ satisfying the following conditions:
\begin{eqnarray}
&&(b b^{'})b^{''}=(b\cd x^1)[(b^{'}\cd x^2)(b^{''}\cd x^3)],
\label{rma1} \\
&&(b b^{'})\cd h=(b\cd h_1)(b^{'}\cd h_2),\;\;\;
1_B\cd h=\va (h)1_B,
\end{eqnarray}
for all $b, b^{'}, b^{''}\in B$ and $h\in H$, where $b\ot h\ra
b\cd h$ is the right $H$-module structure of $B$.\\
${\;\;\;}$
Recall from \cite{hn1} the notion of comodule algebra over a
quasi-bialgebra.
\begin{definition}
Let $H$ be a quasi-bialgebra. A unital associative algebra
$\mathfrak{A}$ is called a right $H$-comodule algebra if there
exist an algebra morphism $\r :\mathfrak{A}\ra \mathfrak{A}\ot H$
and an invertible element $\Phi _{\r }\in \mathfrak{A}\ot H\ot H$
such that:
\renewcommand{\theequation}{\thesection.\arabic{equation}}
\begin{eqnarray}
&&\Phi _{\r }(\r \ot id)(\r (\mf {a}))=(id\ot \Delta
)(\r (\mf {a}))\Phi _{\r }, 
\mbox{${\;\;\;}$$\forall $ $\mf {a}\in
\mathfrak{A}$,}\label{rca1}\\[1mm]%
&&(1_{\mf {A}}\ot \Phi)(id\ot \Delta \ot id)(\Phi _{\r })(\Phi
_{\r }\ot 1_H)= (id\ot id\ot \Delta )(\Phi _{\r })(\r \ot id\ot
id)(\Phi _{\r }),\label{rca2}\\[1mm]%
&&(id\ot \va)\circ \r =id ,\label{rca3}\\[1mm]%
&&(id\ot \va \ot id)(\Phi _{\r })=(id\ot id\ot \va )(\Phi _{\r }
)=1_{\mathfrak{A}}\ot 1_H.\label{rca4}
\end{eqnarray}
Similarly, a unital associative algebra $\mathfrak{B}$ is called
a left $H$-comodule algebra if there exist an algebra morphism $\l
: \mf {B}\ra H\ot \mathfrak{B}$ and an invertible element $\Phi
_{\l }\in H\ot H\ot \mathfrak{B}$ such that:
\begin{eqnarray}
&&(id\ot \l )(\l (\mf {b}))\Phi _{\l }=\Phi _{\l
}(\Delta \ot id)(\l (\mf {b})),
\mbox{${\;\;\;}$$\forall $ $\mf {b}\in \mathfrak{B}$,}
\label{lca1}\\[1mm]%
&&(1_H\ot \Phi _{\l })(id\ot \Delta \ot id)(\Phi _{\l })(\Phi \ot
1_{\mf {B}})= (id\ot id\ot \l )(\Phi _{\l })(\Delta \ot id\ot
id)(\Phi _{\l }),\label{lca2}\\[1mm]%
&&(\va \ot id)\circ \l =id ,\label{lca3}\\[1mm]%
&&(id\ot \va \ot id)(\Phi _{\l })=(\va \ot id\ot id)(\Phi _{\l }
)=1_H\ot 1_{\mathfrak{B}}.\label{lca4}
\end{eqnarray}
\end{definition}
${\;\;\;}$
When $H$ is a quasi-bialgebra, particular examples of left and
right $H$-comodule algebras are given by $\mf {A}=\mf {B}=H$ and
$\r =\l =\Delta $,
$\Phi _{\r }=\Phi _{\l }=\Phi $. Another basic example of a comodule 
algebra is provided by the smash product. Namely, if $A$ is a left 
$H$-module algebra, then $(A\# H, \rho , \Phi _{\rho })$ is a right 
$H$-comodule algebra with structures (cf. \cite{bpv}):
\begin{eqnarray*}
&&\rho :A\# H\rightarrow (A\# H)\otimes H, \;\;\;
\rho (a\# h)=(x^1\cdot a\# x^2h_1)\otimes x^3h_2, \\
&&\Phi _{\rho }=(1\# X^1)\otimes X^2\otimes X^3\in (A\# H)\otimes H
\otimes H.
\end{eqnarray*} 
For a right $H$-comodule algebra $({\mf A}, \r , \Phi _{\r })$ we
will denote, for any $\mfa \in {\mf A}$, by 
$\r (\mfa )=\mfa _{<0>}\ot \mfa _{<1>}$,  $(\r  
\ot id)(\r (\mfa ))=\mfa _{<0, 0>}\ot \mfa _{<0, 1>} \ot \mfa 
_{<1>}$ etc. Similarly, for a left $H$-comodule 
algebra $({\mf B}, \l , \Phi _{\l })$, if $\mfb \in {\mf B}$ then
we denote 
$\l (\mfb )=\mfb _{[-1]}\ot \mfb _{[0]}$,  
$(id\ot \l )(\l (\mfb ))=\mfb _{[-1]}\ot \mfb _{[0,-1]}\ot  
\mfb _{[0, 0]}$ etc. 
In analogy with the notation for the reassociator $\Phi 
$ of $H$, we will write 
\begin{eqnarray*}
&&\Phi _{\r }=\tilde {X}^1_{\r }\ot \tilde {X}^2_{\r }\ot  
\tilde {X}^3_{\r }=
\tilde {Y}^1_{\r }\ot \tilde {Y}^2_{\r }\ot \tilde {Y}^3_{\r }=\cdots  \\
&&\Phi _{\r }^{-1}=\tilde {x}^1_{\r }\ot \tilde {x}^2_{\r }\ot  
\tilde {x}^3_{\r }=\tilde {y}^1_{\r }\ot \tilde {y}^2_{\r 
}\ot \tilde {y}^3_{\r }=\cdots  
\end{eqnarray*}
and similarly for the element $\Phi _{\l }$ of a left $H$-comodule 
algebra $\mf {B}$. \\
${\;\;\;}$
If $\mf {A}$ is a right $H$-comodule algebra then we define the elements 
$\tilde{p}_{\r }, \tilde{q}_{\r }\in {\mf A}\ot H$ as follows:
\begin{equation}\label{tpqr}
\tilde {p}_{\r }=\tilde {p}^1_{\r }\ot \tilde {p}^2_{\r}
=\tx ^1_{\r }\ot \tx ^2_{\r }\b S(\tx ^3_{\r }), 
\mbox{${\;\;\;}$}
\tilde {q}_{\r }=\tilde {q}^1_{\r }\ot \tilde {q}^2_{\r}
=\tX ^1_{\r }\ot \smi (\a \tX ^3_{\r })\tX ^2_{\r }.
\end{equation}
By \cite[Lemma 9.1]{hn1}, we have the following relations, for all $\mfa \in
{\mf A}$:
\begin{eqnarray}
&&
\r (\mfa_{<0>})\tilde {p}_{\r }[1_{\mf A}\ot S(\mfa _{<1>})]
=\tilde {p}_{\r }[\mfa \ot 1_H],\label{tpqr1}\\
&&
[1_{\mf A}\ot \smi (\mfa_{<1>})]\tilde {q}_{\r }\r
(\mfa_{<0>})=[\mfa \ot 1_H]\tilde {q}_{\r },\label{tpqr1a}\\
&&
\r (\tilde {q}^1_{\r })\tilde {p}_{\r }[1_{\mf A}\ot 
S(\tilde{q}^2_{\r })]=1_{\mf A}\ot 1_H\label{tpqr2},\\
&&
[1_{\mf A}\ot \smi (\tilde {p}^2_{\r })]\tilde {q}_{\r }\r
(\tilde {p}^1_{\r })=1_{\mf A}\ot 1_H,\label{tpqr2a}\\
&&
\Phi _{\r }(\r \ot id_H)(\tilde {p}_{\r })(\tilde
{p}_{\r }\ot id_H)\nonumber\\
&&\hspace*{2cm}
=(id_{\mf A}\ot \Delta )
(\r (\tx ^1_{\r })\tilde {p}_{\r  
})(1_{\mf A}\ot g^1S(\tx ^3_{\r })\ot g^2S(\tx ^2_{\r })),\label{relmare} \\
&&
(\tilde {q}_{\r }\ot 1_H)(\r \ot id_H)(\tilde {q}_{\r })
\Phi _{\r }^{-1}\nonumber\\
&&\hspace*{2cm}
=[1_{\mf A}\ot \smi (f^2\tX ^3_{\r })\ot  
\smi (f^1\tX ^2_{\r })](id _{\mf A}\ot \Delta )
(\tilde {q}_{\r }\r (\tX ^1_{\r })),\label{relfoartemare}  
\end{eqnarray}
where $f^{-1}=g^1\otimes g^2$ is given by (\ref{g}). If ${\mf A}$ is 
$H$ itself, the elements 
$\tilde{p}_{\r }, \tilde{q}_{\r }\in H\ot H$ are denoted by $p_R$ and 
$q_R$ and are given by the formulae
\begin{eqnarray}
&&p_R=p^1\ot p^2=x^1\ot x^2\b S(x^3),  
\mbox{${\;\;\;}$}
q_R=q^1\ot q^2=X^1\ot \smi (\a X^3)X^2. \label{psiq}
\end{eqnarray}
The next definition appeared in \cite{hn1} under the 
name "quasi-commuting pair of $H$-coactions".   
\begin{definition}
Let $H$ be a quasi-bialgebra. By an $H$-bicomodule algebra $\mb {A}$ 
we mean a quintuple $(\l, \r , \Phi _{\l }, \Phi _{\r }, \Phi
_{\l , \r })$, where $\l $ and $\r $ are left and right
$H$-coactions on $\mb {A}$, respectively, and where $\Phi _{\l
}\in H\ot H\ot \mb {A}$, $\Phi _{\r }\in \mb {A}\ot H\ot H$ and
$\Phi _{\l , \r }\in H\ot \mb {A}\ot H$ are invertible elements,
such that 
$(\mb {A}, \l , \Phi _{\l })$ is a left $H$-comodule algebra, 
$(\mb {A}, \r , \Phi _{\r })$ is a right $H$-comodule algebra and 
the following compatibility relations hold:
\renewcommand{\theequation}{\thesection.\arabic{equation}}
\begin{eqnarray}
&&\hspace{-1.8cm}\Phi _{\l , \r }(\l \ot id)(\r (u))=(id\ot \r )(\l 
(u))\Phi _{\l, \r }, \mbox{${\;\;}$$\forall $ $u\in \mb  
{A}$,}\label{bca1}\\[1mm]%
&&\hspace{-1.8cm}(1_H\ot \Phi _{\l , \r })(id\ot \l \ot id)(\Phi
_{\l , \r }) (\Phi _{\l }\ot 1_H)=(id\ot id\ot \r )(\Phi _{\l
})(\Delta \ot
id\ot id)(\Phi _{\l , \r }), \label{bca2}\\[1mm]%
&&\hspace{-1.8cm}(1_H\ot \Phi _{\r })(id\ot \r \ot id)(\Phi _{\l ,\r
})(\Phi _{\l , \r }\ot 1_H)= (id\ot id\ot \Delta )(\Phi _{\l , \r
})(\l \ot id\ot id) (\Phi _{\r }).\label{bca3}
\end{eqnarray}
\end{definition}
As pointed out in \cite{hn1}, if $\mb {A}$ is a bicomodule algebra
then, in addition, we have that
\renewcommand{\theequation}{\thesection.\arabic{equation}}
\begin{equation}\label{bca4}
(id_H\ot id_{\mb {A}}\ot \va )(\Phi _{\l , \r })=1_H\ot 1_{\mb
{A}}, \mbox{${\;\;}$} (\va \ot id_{\mb {A}}\ot id_H)(\Phi _{\l ,
\r })= 1_{\mb {A}} \ot 1_H.
\end{equation}
An example of a bicomodule algebra is $\mb {A}=H$, $\l =\r  
=\Delta $ and $\Phi _{\l }=\Phi _{\r }= \Phi _{\l , \r }=\Phi $.
For the left and right comodule algebra structures of $\mb 
{A}$ we will keep the same notation as above. We also denote 
$\Phi _{\l , \r }=\Theta ^1\ot \Theta ^2\ot \Theta ^3=
\ov {\Theta }^1\ot \ov {\Theta }^2\ot \ov {\Theta }^3$ and  
$\Phi ^{-1}_{\l , \r }=\theta ^1\ot \theta ^2\ot \theta
^3=\ov {\theta }^1\ot \ov {\theta }^2\ot \ov {\theta }^3$.\\
${\;\;\;}$ 
Let us denote by $_H{\cal M}_H$ the category of $H$-bimodules; it is 
also a monoidal category, the 
associativity constraints being given by ${\bf a^{'}}_{U, V, W}:
(U\ot V)\ot W\ra U\ot (V\ot W)$,  
${\bf a^{'}}_{U, V, W}((u\ot v)\ot w)= \Phi \cd (u\ot (v\ot w))\cd
\Phi ^{-1}$,  
for $U, V, W\in {}_H{\cal M}_H$ and $u\in U$, $v\in V$, $w\in 
W$. Therefore, we can define algebras in the category of
$H$-bimodules, which will be called $H$-bimodule 
algebras. Namely, a $k$-vector space ${\cal A}$ is an 
$H$-bimodule algebra if ${\cal A}$ is an $H$-bimodule (denote
the actions by $h\cd \varphi $ and $\varphi \cd h$, for $h\in H$ and 
$\varphi \in {\cal A}$)   
with a multiplication and a usual unit $1_{\cal A}$ such that:
\renewcommand{\theequation}{\thesection.\arabic{equation}}
\begin{eqnarray}
&&\hspace{-2cm}(\varphi \psi )\xi =(X^1\cd \varphi \cd
x^1)[(X^2\cd \psi \cd x^2)(X^3\cd \xi \cd x^3)],
\mbox{${\;\;}$$\forall $ $\varphi , \psi , \xi \in
{\cal A}$,}\label{bma1}\\[1mm]%
&&\hspace{-2cm}h\cd (\varphi \psi)=(h_1\cd \varphi)(h_2\cd
\psi ), 
\mbox{${\;\;}$} 
(\varphi \psi )\cd h=(\varphi \cd h_1)(\psi \cd h_2),  
\mbox{${\;\;}$$\forall $ $\varphi , \psi \in {\cal A}$, $h\in H$,} 
\label{bma2}\\
&&\hspace{-2cm}
h\cd 1_{\cal A}=\va (h)1_{\cal A},
\mbox{${\;\;}$}
1_{\cal A}\cd h=\va (h)1_{\cal A},
\mbox{${\;\;}$$\forall \; h\in H$.}\label{bma3}
\end{eqnarray}
If $H$ is a quasi-bialgebra then $H^*$,  
the linear dual of $H$, is an $H$-bimodule with  
$H$-actions 
$<h\rh \v , h'>=\v (h'h)$ and 
$<\v \lh h, h'>=\v (hh')$, for all $\varphi \in H^*$, $h, h'\in H$. 
The convolution $<\v \psi , h>=\sum \v (h_1)\psi (h_2)$, for $\v , 
\psi \in H^*$ and $h\in H$, is a multiplication on $H^*$, 
and with this multiplication $H^*$ becomes an $H$-bimodule algebra. 
Note that a left (right) $H$-module algebra becomes also an $H$-bimodule 
algebra, with trivial right (left) $H$-action. Unlike the case of 
ordinary bialgebras, an $H$-bimodule algebra is {\it not} necessarily 
a left or right $H$-module algebra, for instance there is no visible left 
or right $H$-module algebra structure on $H^*$.\\
${\;\;\;}$
We recall from \cite{majid} the definition of (left) Yetter-Drinfeld modules 
over a quasi-bialgebra $H$. 
\begin{definition}
A $k$-linear space $M$ is called a left Yetter-Drinfeld module over $H$ if 
$M$ is a left $H$-module (with action denoted by $h\ot m\mapsto h\cdot m$) 
and $H$ coacts on $M$ to the left (the coaction is denoted by 
$\lambda _M:M\rightarrow H\ot M$, $\lambda _M(m)=m_{(-1)}\ot m_{(0)}$) such 
that
\begin{eqnarray}
&&X^1m_{(-1)}\ot (X^2\cdot m_{(0)})_{(-1)}X^3\ot (X^2\cdot m_{(0)})_{(0)}
\nonumber \\
&&\;\;\;\;\;\;\;\;\;
=X^1(Y^1\cdot m)_{(-1)_1}Y^2\ot X^2(Y^1\cdot m)_{(-1)_2}Y^3\ot  
X^3\cdot (Y^1\cdot m)_{(0)}, \label{yd1} \\
&&\varepsilon (m_{(-1)})m_{(0)}=m, \label{yd2} \\
&&h_1m_{(-1)}\ot h_2\cdot m_{(0)}=(h_1\cdot m)_{(-1)}h_2\ot 
(h_1\cdot m)_{(0)}, \label{yd3}
\end{eqnarray}
for all $m\in M$ and $h\in H$. 
The category $_H^H{\cal YD}$ consists of such objects, the morphisms in the 
category being the $H$-linear maps intertwining the $H$-coactions. 
\end{definition}
${\;\;\;}$
The category $_H^H{\cal YD}$ is (pre) braided monoidal; explicitly, 
if $(M, \lambda _M)$ and $(N, \lambda _N)$ are objects in $_H^H{\cal YD}$, 
then $(M\ot N, \lambda _{M\ot N})$ is also object in 
$_H^H{\cal YD}$, where $M\ot N$ is a left $H$-module with action $h\cdot 
(m\ot n)=h_1\cdot m\ot h_2\cdot n$, and the coaction $\lambda _{M\ot N}$ is 
given by 
\begin{eqnarray*}
&&\lambda _{M\ot N}(m\ot n)=X^1(x^1Y^1\cdot m)_{(-1)}x^2(Y^2\cdot n)_{(-1)}
Y^3\ot X^2\cdot (x^1Y^1\cdot m)_{(0)}\\ 
&&\;\;\;\;\;\;\;\;\;\;\;\;\;\;\;\;\;\;\;\;\;\;\;\;\;\;\;\;\;\;\;\;
\ot X^3x^3\cdot (Y^2\cdot n)_{(0)}.
\end{eqnarray*}
The associativity constraints are the same as in $_H{\cal M}$, and the 
(pre) braiding is given by 
\begin{eqnarray*}
&&c_{M, N}:M\ot N\rightarrow N\ot M, \;\;c_{M, N}(m\ot n)=m_{(-1)}\cdot n\ot 
m_{(0)}.
\end{eqnarray*}
${\;\;\;}$ 
Since $_H^H{\cal YD}$ is a monoidal category, we can speak about algebras 
in $_H^H{\cal YD}$. Namely, if $A$ is an object in $_H^H{\cal YD}$, then 
$A$ is an algebra in $_H^H{\cal YD}$ if and only if 
$A$ is a left $H$-module algebra and 
$A$ is a left quasi-comodule algebra, that is its unit and 
multiplication intertwine the $H$-coaction $\lambda _A$, namely 
(for all $a, a'\in A$):  
\begin{eqnarray}
&&\lambda _A(1_A)=1_H\ot 1_A, \label{unitate} \\
&&\lambda _A(aa')=X^1(x^1Y^1\cdot a)_{(-1)}x^2(Y^2\cdot a')_{(-1)}Y^3 
\nonumber \\
&&\;\;\;\;\;\;\;\;\;\;\;\;\;\;\;\;\;\;\;\;
\ot [X^2\cdot (x^1Y^1\cdot a)_{(0)}][X^3x^3\cdot (Y^2\cdot a')_{(0)}].
\label{multi}
\end{eqnarray}
${\;\;\;}$
If $H$ is a quasi-Hopf algebra, we can consider the algebra map 
$id_H:H\rightarrow H$ and then the left $H$-module algebra $H^{id}$, which 
was denoted by $H_0$ in \cite{bpv}; its unit is $\beta $, the $H$-action 
is $h\triangleright h'=h_1h'S(h_2)$, and the multiplication is 
$h*h'=X^1hS(x^1X^2)\alpha x^2X^3_1h'S(x^3X^3_2)$. 
Moreover, it was proved in \cite{bn} that $H_0$ becomes an algebra  
in $_H^H{\cal YD}$, with coaction $\lambda _{H_0}:H_0\rightarrow H\otimes 
H_0$, 
\begin{eqnarray*}
&&\lambda _{H_0}(h):=h_{(-1)}\otimes h_{(0)}=X^1Y^1_1h_1g^1S(q^2Y^2_2)Y^3
\otimes X^2Y^1_2h_2g^2S(X^3q^1Y^2_1), 
\end{eqnarray*}
where $f^{-1}=g^1\otimes g^2$ is given by (\ref{g}) and  
$q_R=q^1\otimes q^2$ is given by (\ref{psiq}).\\
${\;\;\;}$
Let $H$ be a quasi-bialgebra, $A$ a left $H$-module algebra and
$\mathfrak{B}$ a left $H$-comodule algebra. Denote by 
$A\gsm \mathfrak{B}$ the $k$-vector space $A\ot\mathfrak{B}$ with 
multiplication:
\begin{eqnarray*}
(a\gsm \mf {b})(a'\gsm \mf {b}')=(\tilde {x}^1_{\l }\cd
a)(\tilde {x}^2_{\l }\mf {b}_{[-1]}\cd a')\gsm 
\tilde {x}^3_{\l }\mf {b}_{[0]}\mf {b}' ,\;\;\;
\forall \;a, a'\in A, \; \mf {b}, \mf {b}'\in \mathfrak{B}.
\end{eqnarray*}
By \cite{bc}, $A\gsm \mathfrak{B}$ is an associative algebra 
with unit $1_A\gsm 1_{\mathfrak{B}}$, called the generalized smash product 
of $A$ and $\mathfrak{B}$ (it coincides with the usual smash product if 
$\mathfrak{B}=H$). \\
${\;\;\;}$ 
We recall the so-called L-R-smash product, introduced in \cite{pv} as 
a generalization of the cocommutative case from \cite{bbm}, \cite{bggs}.  
Let $H$ be a quasi-bialgebra, ${\cal A}$ an $H$-bimodule 
algebra and ${\mb A}$ an $H$-bicomodule algebra. Define on ${\cal A} 
\ot {\mb A}$ the product  
\begin{eqnarray}\label{nat} 
&&(\varphi \nat u)(\varphi ' \nat u')=
(\tilde{x}^1_{\lambda }\cd \varphi \cd 
\theta ^3u'_{<1>}\tilde{x}^2_{\rho })
(\tilde{x}^2_{\lambda }u_{[-1]}\theta ^1\cd \varphi ' \cd 
\tilde{x}^3_{\rho })\nat \tilde{x}^3_{\lambda }u_{[0]}\theta ^2u'_{<0>}
\tilde{x}^1_{\rho }
\end{eqnarray}
for $\v , \v ' \in {\cal A}$ and $u, u'\in {\mb A}$,    
where we write $\v \nat u$ in place of $\v \ot u$      
to distinguish the new algebraic structure.    
Then this product defines on 
${\cal A}\ot {\mb A}$ a structure of associative algebra with unit 
$1_{\cal A}\nat 1_{\mb A}$, denoted by     
${\cal A}\nat {\mb A}$ and called the L-R-smash product. Note that if 
the right $H$-action on $\cal A$ is trivial then 
${\cal A}\nat {\mb A}$ coincides with the generalized smash product 
${\cal A}\gsm {\mb A}$.\\
${\;\;\;}$
If H is a quasi-Hopf algebra, ${\cal A}$ an $H$-bimodule   
algebra and ${\mb A}$ an $H$-bicomodule algebra, there exists another 
associative algebra structure built on ${\cal A}\otimes {\mb A}$, which was 
introduced in \cite{bpvo} under the name 
generalized diagonal crossed product (for ${\cal A}=H^*$ it gives  
the diagonal crossed product introduced by Hausser and Nill in \cite{hn1}). 
It is constructed as follows. Define the element   
\begin{eqnarray}
&&
\hspace*{-1cm}
\O =(\tX ^1_{\r })_{[-1]_1}\tx ^1_{\l }\theta ^1\ot (\tX
^1_{\r })_{[-1]_2}\tx ^2_{\l }\theta ^2_{[-1]}\ot (\tX ^1_{\r
})_{[0]}\tx ^3_{\l }\theta ^2_{[0]}\ot \smi (f^1\tX ^2_{\r }\theta
^3)\ot \smi (f^2\tX ^3_{\r }) \label{formomega}
\end{eqnarray}
in $H^{\ot 2}\ot {\mb A}\ot H^{\ot 2}$, 
where $\Phi _{\r }=\tX ^1_{\r }\ot \tX ^2_{\r }\ot \tX ^3_{\r 
}$, $\Phi _{\l }^{-1}=\tx ^1_{\l }\ot \tx ^2_{\l }\ot \tx 
^3_{\l }$,  
$\Phi _{\l , \r }^{-1}=\theta ^1\ot \theta ^2\ot  
\theta ^3$ and $f=f^1\ot f^2$ is the twist defined in (\ref{f}). Then  
define a multiplication on ${\cal A}\ot {\mb A}$, by: 
\begin{eqnarray}\label{gdp}
(\varphi \bowtie u)(\varphi ' \bowtie u')=(\O ^1\cd \varphi \cd    
\O ^5)(\O ^2u_{<0>_{[-1]}}\cd \varphi ' \cd \smi (u_{<1>})\O ^4)\bowtie 
\O ^3u_{<0>_{[0]}}u', 
\end{eqnarray}
for all $\varphi , \varphi ' \in {\cal A}$ and $u, u'\in {\mb A}$, where 
we denoted $\Omega =\Omega ^1\otimes ... \otimes \Omega ^5$. Then this   
multiplication defines an associative algebra structure with unit  
$1_{\cal A}\bowtie 1_{\mb A}$, which will be denoted by 
${\cal A}\bowtie {\mb A}$. Note that, as in the case of the L-R-smash 
product, if 
the right $H$-action on $\cal A$ is trivial then 
${\cal A}\bowtie {\mb A}$ coincides with the generalized smash product  
${\cal A}\gsm {\mb A}$.\\ 
${\;\;\;}$It was proved in \cite{pv} that actually ${\cal A}\nat {\mb A}$ 
and ${\cal A}\bowtie {\mb A}$ are isomorphic as algebras, a pair of 
algebra isomorphisms $\nu :{\cal A}\bowtie {\mb A}\rightarrow 
{\cal A}\nat {\mb A}$ and 
$\nu ^{-1}:{\cal A}\nat {\mb A}\rightarrow {\cal A}\bowtie {\mb A}$ being 
given by  
\begin{eqnarray}
&&\nu (\varphi \bowtie u)=\Theta ^1\cdot \varphi \cdot 
\smi (\Theta ^3)\tqrb \Theta ^2_{<1>}u_{<1>}\nat 
\tqra \Theta ^2_{<0>}u_{<0>}, \label{formulanu} \\
&&\nu ^{-1}(\varphi \nat u)=\theta ^1\cdot \varphi \cdot 
S^{-1}(\theta ^3u_{<1>}\tilde{p}^2_{\rho })\bowtie \theta ^2u_{<0>}
\tilde{p}^1_{\rho }, \label{formulanuinvers}
\end{eqnarray}
for all $\varphi \in {\cal A}$, $u\in {\mb A}$, where 
$\tilde {p}_{\r }=\tilde {p}^1_{\r }\ot \tilde {p}^2_{\r}$ and  
$\tilde {q}_{\r }=\tilde {q}^1_{\r }\ot \tilde {q}^2_{\r}$ are the elements 
given by (\ref{tpqr}).\\
${\;\;\;}$
We recall now several facts about twisted tensor products of algebras 
in monoidal categories (see for instance \cite{cap}, \cite{vandaele}). 
Let ${\cal C}$ be a monoidal category with  
associativity constraints $a_{U, V, W}: (U\ot V)\ot W\rightarrow 
U\ot (V\ot W)$ and unit $I$. If $A$ is an algebra in ${\cal C}$ we 
denote its multiplication by $\mu _A:A\ot A\rightarrow A$ and its unit 
by $\eta _A:I\rightarrow A$. Let $(A, \mu _A, \eta _A)$ and 
$(B, \mu _B, \eta _B)$ be two algebras in ${\cal C}$ and 
$R:B\ot A\rightarrow A\ot B$ a morphism in ${\cal C}$. We call $R$ a 
{\it twisting map} between $A$ and $B$ if  
the following conditions hold: 
\begin{eqnarray}
&&R\circ (id_B\ot \eta _A)=\eta _A\ot id_B,\;\;\;
R\circ (\eta _B\ot id_A)=id_A\ot \eta _B, \label{unitw} \\
&&R\circ (\mu _B\ot id_A)=(id_A\ot \mu _B)\circ a_{A, B, B}\circ  
(R\ot id_B)\circ a_{B, A, B}^{-1}\nonumber \\
&&\hspace*{4cm}
\circ (id_B\ot R)\circ a_{B, B, A}, 
\label{twist1} \\
&&R\circ (id_B\ot \mu _A)=(\mu _A\ot id_B)\circ a_{A, A, B}^{-1}\circ  
(id_A\ot R)\circ a_{A, B, A}\nonumber \\
&&\hspace*{4cm}
\circ (R\ot id_A)\circ a_{B, A, A}^{-1}. 
\label{twist2}
\end{eqnarray}   
Given such a twisting map, $A\ot B$ becomes an  
algebra in ${\cal C}$, with multiplication 
\begin{eqnarray}
\mu &=&(\mu _A\ot \mu _B)\circ a_{A, A, B\ot B}^{-1}\circ 
(id_A\ot a_{A, B, B})\circ (id_A\ot R\ot id_B)\nonumber \\
&&\circ (id_A\ot a_{B, A, B}^{-1})\circ a_{A, B, A\ot B} \label{multitwisting}
\end{eqnarray}
and unit $\eta =\eta _A\ot \eta _B$. This algebra structure on $A\ot B$ 
will be denoted by $A\ot _RB$ and will be called the {\it twisted tensor  
product} of $A$ and $B$ afforded by the twisting map $R$. It has 
moreover the property that the morphisms $i_A:=id_A\ot \eta _B:A\rightarrow 
A\ot _RB$ and $i_B:=\eta _A\ot id_B:B\rightarrow A\ot _RB$ are 
morphisms of algebras in ${\cal C}$.\\
${\;\;\;}$
If $c_{U, V}:U\ot V\rightarrow V\ot U$ is a braiding on ${\cal C}$, then  
for any two algebras $A$ and $B$ in ${\cal C}$, the map 
$R=c_{B, A}:B\ot A\rightarrow A\ot B$ is a twisting map, and in this case 
$A\ot _RB$ is denoted by $A\underline{\ot}B$ and 
is called the {\it braided tensor product} of $A$ and $B$ (relative to the 
braiding $c$).\\
${\;\;\;}$
If $A\ot _RB$ is a twisted tensor product of algebras in ${\cal C}$, 
it has the following universal property (see \cite{bd2}, \cite{cimz}): if  
$(X, \mu _X, \eta _X)$ is an algebra in ${\cal C}$ and $u:A\rightarrow X$ 
and $v:B\rightarrow X$ are morphisms of algebras in ${\cal C}$ such that 
\begin{eqnarray}
&&\mu _X\circ (u\ot v)\circ R=\mu _X\circ (v\ot u), \label{comgen}
\end{eqnarray}
then there exists a unique morphism $w:A\ot _RB\rightarrow X$ of algebras 
in ${\cal C}$ such that $w\circ i_A=u$ and $w\circ i_B=v$. This 
morphism $w$ is given explicitely by $w=\mu _X\circ (u\ot v)$.   
\section{Twisted tensor products of  
quasialgebras}\label{sec2}
\setcounter{equation}{0}
\begin{proposition} \label{ambeleparti}
Let $H$ be a quasi-bialgebra, ${\cal A}$ an $H$-bimodule algebra and 
$A$ an algebra in $_H^H{\cal YD}$. We regard $A$ as an $H$-bimodule 
algebra with trivial right $H$-action. Define the map 
\begin{eqnarray}
&&R:A\ot {\cal A}\rightarrow {\cal A}\ot A, \;\;\;\; 
R(a\ot \varphi )=a_{(-1)}\cdot \varphi \ot a_{(0)}, \;\;\; \forall \;  
a\in A,\;\varphi \in {\cal A}. \label{Rgen}
\end{eqnarray}
Then $R$ is a twisting map between ${\cal A}$ and $A$ in the monoidal 
category $_H{\cal M}_H$. We will denote by ${\cal A}\odot A$ the $H$-bimodule 
algebra ${\cal A}\ot _RA$, the twisted tensor product of ${\cal A}$ and 
$A$ in $_H{\cal M}_H$.
\end{proposition} 
\begin{proof}
The fact that $R$ is right $H$-linear is obvious, and the fact that it is left 
$H$-linear follows immediately from (\ref{yd3}), so indeed $R$ is a 
morphism in $_H{\cal M}_H$. The relations (\ref{unitw}) follow immediately 
from (\ref{yd2}) and (\ref{unitate}).  
The relation (\ref{twist1}) reduces to (\ref{multi}), while the relation 
(\ref{twist2}) reduces to (\ref{yd1}).
\end{proof}
\begin{remark}
The explicit structure of ${\cal A}\odot A$ is the following: the 
unit is $1_{{\cal A}}\ot 1_A$, the left $H$-action is $h\cdot (\varphi \ot a)
=h_1\cdot \varphi \ot h_2\cdot a$, the right $H$-action is 
$(\varphi \ot a)\cdot h=\varphi \cdot h\ot a$, and the multiplication is 
(by (\ref{multitwisting})):
\begin{eqnarray*}
(\varphi \ot a)(\varphi '\ot a')&=&(y^1X^1\cdot \varphi )(y^2Y^1
(x^1X^2\cdot a)_{(-1)}x^2X^3_1\cdot \varphi ')\\
&&\ot (y^3_1Y^2\cdot (x^1X^2\cdot a)_{(0)})(y^3_2Y^3x^3X^3_2\cdot a').
\end{eqnarray*}
\end{remark}
${\;\;\;}$
We will be mainly interested in the following particular case of 
Proposition \ref{ambeleparti}. 
\begin{corollary}
Let $H$ be a quasi-Hopf algebra, $A$ a left $H$-module algebra and 
$H_0$ the algebra in $_H^H{\cal YD}$ as in the Preliminaries. Define the 
map $R:H_0\ot A\rightarrow A\ot H_0$, 
\begin{eqnarray}
&&R(h\ot a)=h_{(-1)}\cdot a\ot h_{(0)}=X^1Y^1_1h_1g^1S(q^2Y^2_2)Y^3\cdot a
\ot X^2Y^1_2h_2g^2S(X^3q^1Y^2_1). \label{Rsmash}
\end{eqnarray}
Then $R$ is a twisting map between $A$ and $H_0$ in the monoidal category 
$_H{\cal M}$. We will denote by $A\diamond H_0$ the left $H$-module algebra 
$A\ot _RH_0$, the twisted tensor product of $A$ and $H_0$ in $_H{\cal M}$. 
Its unit is $1_A\ot \beta $, the $H$-action is $h\cdot (a\ot h')=
h_1\cdot a\ot h_2\triangleright h'$, and the multiplication is
\begin{eqnarray*}
(a\ot h)(a'\ot h')&=&(y^1X^1\cdot a)(y^2Y^1
(x^1X^2\triangleright h)_{(-1)}x^2X^3_1\cdot a')\\
&&\ot (y^3_1Y^2\triangleright (x^1X^2\triangleright h)_{(0)})
*(y^3_2Y^3x^3X^3_2\triangleright h').
\end{eqnarray*}
\end{corollary}
${\;\;\;}$
More generally, if $C$ is a left $H$-module algebra and $A$ is an algebra 
in $_H^H{\cal YD}$, then $C\odot A$ is a left $H$-module algebra, which will 
be denoted by $C\diamond A$.\\
${\;\;\;}$  
Since the braiding in $_H^H{\cal YD}$ is given by 
$m\ot n\mapsto m_{(-1)}\cdot n\ot m_{(0)}$, we obtain:
\begin{corollary}
Let $H$ be a quasi-Hopf algebra and $A$ an algebra in $_H^H{\cal YD}$. 
Then the left $H$-module algebra $A\diamond H_0$ is actually an algebra 
in $_H^H{\cal YD}$ and it coincides with the braided tensor product 
$A\underline{\ot }H_0$ in $_H^H{\cal YD}$.
\end{corollary}
${\;\;\;}$As another class of examples of twisted tensor products of 
quasialgebras, we present the quasialgebra analogue of the 
so-called {\it Clifford process}, cf. \cite{am2}, \cite{wene}. Let $A$ be a 
(not necessarily associative) algebra over $k$ with unit $1$, let $q\in k$ 
be a nonzero fixed element and $\sigma :A\rightarrow A$ an involutive 
unital algebra automorphism. We denote by $C(k, q)$ the $2$-dimensional 
associative algebra $k[v]/(v^2=q)$. We define the linear map 
\begin{eqnarray}
&&R:C(k, q)\otimes A\rightarrow A\otimes C(k, q), \;\;\;
R(1\otimes a)=a\otimes 1, \;\;R(v\otimes a)=\sigma (a)\otimes v, 
\label{RClifford}
\end{eqnarray}
for all $a\in A$. The Clifford process associates to the pair  
$(A, \sigma )$ a  
(not necessarily associative) algebra structure built on $A\otimes C(k, q)$, 
with multiplication 
\begin{eqnarray}
&&(a\otimes 1+b\otimes v)(c\otimes 1+d\otimes v)=(ac+qb\sigma (d))\otimes 1
+(ad+b\sigma (c))\otimes v,  
\label{multClifford}
\end{eqnarray}
for all $a, b, c, d\in A$. This algebra structure is denoted by 
$\overline{A}$. Moreover, the linear map  
\begin{eqnarray*}
&&\overline{\sigma }:\overline{A}\rightarrow \overline{A}, \;\;\;
\overline{\sigma }(a\otimes 1+b\otimes v)=\sigma (a)\otimes 1-\sigma (b)
\otimes v, \;\;\;\forall \;a, b\in A,
\end{eqnarray*}
is an involutive unital algebra automorphism for $\overline{A}$. It is clear 
(cf. \cite{am2}) that if $A$ is associative then $\overline{A}$ is also 
associative, as in this case $R$ is a twisting map and $\overline{A}$ is the 
twisted tensor product $\overline{A}=A\otimes _RC(k, q)$. \\
${\;\;\;}$Consider now the initial data for a Clifford process, and 
assume that $A$ is a left module algebra over a quasi-bialgebra $H$  
and $\sigma $ is $H$-linear. We regard $C(k, q)$ as a left 
$H$-module algebra with trivial $H$-action. Then one can see that the map 
$R$ given by (\ref{RClifford}) is a twisting map between $A$ and 
$C(k, q)$ in the monoidal category $_H{\cal M}$, and so 
$\overline{A}=A\otimes _RC(k, q)$ is also a left $H$-module algebra. It 
is also clear that the extended involutive automorphism 
$\overline{\sigma }$ is $H$-linear too. \\
${\;\;\;}$Dually, a similar result holds if $A$ is a comodule algebra over a 
coquasi-bialgebra. In particular, if $G$ is a group and $A$ is a 
$G$-graded quasialgebra as in \cite{am1} and $\sigma $ is a graded involutive 
automorphism of $A$, then $\overline{A}$ is also a $G$-graded quasialgebra. \\
${\;\;\;}$We emphasize an important conceptual  
difference between the associative and   
quasiassociative versions of the Clifford process. In the associative case 
(and assuming $char (k)\neq 2$), as noted in \cite{am1}, the algebra $A$ 
becomes a $\mathbb{Z}_2$-graded algebra and the  
twisted tensor product  
$A\otimes _RC(k, q)$ is actually a braided tensor product, i. e. the  
twisting map $R$ is obtained from the canonical braiding 
$c(x\otimes y)=(-1)^{|x||y|}y\otimes x$ of the category of 
$\mathbb{Z}_2$-graded  
vector spaces. In the quasiassociative case, the twisting map $R$ 
given by (\ref{RClifford}) does {\it not} seem to come from any 
braiding.  
\section{An isomorphism $A\diamond H_0\simeq (A\# H)^j$ of algebras 
in $_H{\cal M}$}\label{sec3}
\setcounter{equation}{0}
\begin{lemma} \label{comut}
Let $H$ be a quasi-Hopf algebra and $A$ a left $H$-module algebra. Then 
for the maps $j$ and $i_0$ given by (\ref{j}) and (\ref{i0}) the 
following relation holds, for all $h\in H$ and $a\in A$:
\begin{eqnarray}
&&j(h)*i_0(a)=i_0(Y^1X^1_1h_1g^1S(T^2X^2_2)\alpha T^3X^3\cdot a)*
j(Y^2X^1_2h_2g^2S(Y^3T^1X^2_1)), \label{comutare}
\end{eqnarray}
where $f^{-1}=g^1\otimes g^2$ is given by 
(\ref{g}) and $*$ is the multiplication in $(A\# H)^j$. 
\end{lemma}  
\begin{proof}
We compute:
\begin{eqnarray*}
j(h)*i_0(a)&=&(1\# h)*(x^1\cdot a\# x^2\beta S(x^3))\\
&=&(1\# X^1)(1\# h)(1\# S(y^1X^2)\alpha y^2X^3_1)\\
&&(x^1\cdot a\# x^2\beta S(x^3))(1\# S(y^3X^3_2))\\
&=&X^1_1h_1S(y^1X^2)_1\alpha _1y^2_1X^3_{(1, 1)}x^1\cdot a\\
&&\# X^1_2h_2S(y^1X^2)_2\alpha _2y^2_2X^3_{(1, 2)}x^2\beta S(y^3X^3_2x^3)\\
&{{\rm (\ref{q1}, \;\ref{moref}, \;\ref{ca})}\atop =}&
X^1_1h_1g^1S(y^1_2X^2_2)\gamma ^1y^2_1x^1X^3_1\cdot a\\
&&\# X^1_2h_2g^2S(y^1_1X^2_1)\gamma ^2y^2_2x^2X^3_{(2, 1)}\beta  
S(X^3_{(2, 2)})S(y^3x^3)\\
&{{\rm (\ref{q5}, \;\ref{gd})}\atop =}&
X^1_1h_1g^1S(T^2t^1_2y^1_2X^2_2)\alpha T^3t^2y^2_1x^1X^3\cdot a\\
&&\# X^1_2h_2g^2S(T^1t^1_1y^1_1X^2_1)\alpha t^3y^2_2x^2\beta S(y^3x^3)\\
&{{\rm (\ref{q3}, \;\ref{q1}, \;\ref{q5})}\atop =}&
X^1_1h_1g^1S(T^2X^2_2)\alpha T^3X^3\cdot a
\# X^1_2h_2g^2S(T^1X^2_1)S(y^1)\alpha y^2\beta S(y^3)\\
&{{\rm (\ref{q6})}\atop =}&
X^1_1h_1g^1S(T^2X^2_2)\alpha T^3X^3\cdot a
\# X^1_2h_2g^2S(T^1X^2_1)\\
&{{\rm (\ref{inainte})}\atop =}&
i_0(Y^1X^1_1h_1g^1S(T^2X^2_2)\alpha T^3X^3\cdot a)*
j(Y^2X^1_2h_2g^2S(Y^3T^1X^2_1)), 
\end{eqnarray*}
finishing the proof.
\end{proof}
\begin{proposition} \label{izomleft}
Let $H$ be a quasi-Hopf algebra and $A$ a left $H$-module algebra. Then the 
map $\Pi :A\diamond H_0\rightarrow (A\# H)^j$, 
$\;\Pi (a\ot h)=x^1\cdot a\# x^2hS(x^3)$, for all $\;a\in A$ and $h\in H$,   
is an isomorphism of left $H$-module algebras. 
\end{proposition}  
\begin{proof}
Note first that $\Pi $ is bijective, with inverse given by 
$\Pi ^{-1}(a\# h)=X^1\cdot a\ot X^2hS(X^3)$. For proving that $\Pi $ 
is a morphism of left $H$-module algebras, we will use the universal 
property of $A\diamond H_0$ as a twisted tensor product of algebras 
in $_H{\cal M}$. Namely, we know that $i_0:A\rightarrow (A\# H)^j$ is a 
morphism of left $H$-module algebras, and it is easy to see that 
$j:H_0\rightarrow (A\# H)^j$ is also a morphism of left $H$-module algebras.  
Moreover, the relation (\ref{comgen}) reduces in this case exactly to 
(\ref{comutare}). We can thus use the universal property of 
$A\diamond H_0$, which provides a morphism $w:A\diamond H_0\rightarrow 
(A\# H)^j$ of left $H$-module algebras, which is moreover given by 
$w(a\ot h)=i_0(a)*j(h)$. Relation (\ref{inainte}) shows that actually we 
have $w=\Pi $, finishing the proof. 
\end{proof}
${\;\;\;}$
As an application of Proposition \ref{izomleft}, we will prove a certain 
kind of invariance under Drinfeld twisting for the $H$-module algebra 
$A\diamond H_0$. We recall first the following results from \cite{pvo}:
\begin{lemma} (\cite{pvo}) \label{lemaimpor}
Let $H$ be a quasi-Hopf algebra, $B$, $C$ associative algebras,  
$\eta :B\rightarrow C$, $j:H\rightarrow B$, $v:H\rightarrow C$ algebra maps 
such that $\eta \circ j=v$. Then the map $\eta :B^j\rightarrow C^v$  is a 
morphism of left $H$-module algebras. 
\end{lemma}
\begin{proposition} (\cite{pvo}) \label{propajut}
Let $H$ be a quasi-Hopf  
algebra, $F\in H\otimes H$ a gauge transformation, $B$ an associative 
algebra, $v:H\rightarrow B$ an algebra map, which will be denoted by $v_F$ 
when considered as a map from $H_F$ to $B$. Then the map 
\begin{eqnarray*}
&&\psi :(B^v)_{F^{-1}}\rightarrow B^{v_F}, \;\;
\psi (b)=v(F^1)bv(S(F^2)), \;\;\forall \;\;b\in B, 
\end{eqnarray*}
is an isomorphism of left $H_F$-module algebras. 
\end{proposition}
${\;\;\;}$
We can state now the desired result.
\begin{proposition}
Let $H$ be a quasi-Hopf algebra, $F\in H\otimes H$ a gauge transformation and 
$A$ a left $H$-module algebra. Then $(A\diamond H_0)_{F^{-1}}\simeq 
A_{F^{-1}}\diamond (H_F)_0$ as left $H_F$-module algebras. 
\end{proposition}
\begin{proof}
We denote as before by $j:H\rightarrow A\# H$, $j(h)=1\# h$, and by $j_F$ the 
same map when considered as a map from $H_F$ to $A\# H$ or to 
$A_{F^{-1}}\# H_F$. Since obviously the algebra isomorphism 
$\pi :A\# H\simeq A_{F^{-1}}\# H_F$, $\pi (a\# h)=F^1\cdot a\# F^2h$, 
satisfies $\pi \circ j_F=j_F$, we obtain that $(A\# H)^{j_F}\simeq 
(A_{F^{-1}}\# H_F)^{j_F}$ as left $H_F$-module algebras, by using 
Lemma \ref{lemaimpor}. On the other hand, by Proposition \ref{propajut}, 
we obtain that $((A\# H)^j)_{F^{-1}}\simeq (A\# H)^{j_F}$ as left 
$H_F$-module algebras. Now by Proposition \ref{izomleft} we have 
$A\diamond H_0\simeq (A\# H)^j$ as left $H$-module algebras and 
$A_{F^{-1}}\diamond (H_F)_0\simeq (A_{F^{-1}}\# H_F)^{j_F}$ as left 
$H_F$-module algebras. It is also easy to see that if $B$ and $C$ are 
isomorphic left $H$-module algebras then $B_{F^{-1}}$ and 
$C_{F^{-1}}$ are isomorphic left $H_F$-module algebras. Using all these 
facts we finally obtain that $(A\diamond H_0)_{F^{-1}}$ and  
$A_{F^{-1}}\diamond (H_F)_0$ are isomorphic as left $H_F$-module algebras.  
\end{proof}
\section{An isomorphism $A\underline{\otimes}H_0\simeq (A\# H)^j$ of 
algebras in $_H^H{\cal YD}$}\label{sec4}
\setcounter{equation}{0}
${\;\;\;}$
We have seen in the previous section that, if $H$ is a quasi-Hopf algebra and 
$A$ is a left $H$-module algebra, then the map 
$\Pi :A\diamond H_0\rightarrow (A\# H)^j$, 
$\;\Pi (a\ot h)=x^1\cdot a\# x^2hS(x^3)$, 
is an isomorphism of left $H$-module algebras. If we assume that $A$ is 
moreover an algebra in $_H^H{\cal YD}$, in which case $A\diamond H_0$ becomes 
the braided tensor product $A\underline{\ot }H_0$, an algebra in 
$_H^H{\cal YD}$, then we intend to show that $(A\# H)^j$ becomes also 
an algebra in $_H^H{\cal YD}$ in a natural way and $\Pi $ becomes an 
isomorphism of algebras in $_H^H{\cal YD}$. 
We begin with a result of independent interest.
\begin{proposition}\label{smlca}
Let $H$ be a quasi-bialgebra and $A$ an algebra in $_H^H{\cal YD}$. Then 
$(A\# H, \lambda , \Phi _{\lambda })$ is a left $H$-comodule algebra, 
with structures:
\begin{eqnarray*}
&&\lambda :A\# H\rightarrow H\ot (A\# H), \;\;\;\;
\lambda (a\# h)=T^1(t^1\cdot a)_{(-1)}t^2h_1\ot (T^2\cdot 
(t^1\cdot a)_{(0)}\# T^3t^3h_2),\\
&&\Phi _{\lambda }=X^1\ot X^2\ot (1\# X^3)\in H\ot H\ot (A\# H).
\end{eqnarray*}
\end{proposition}
\begin{proof}
We first check that $\lambda $ is an algebra map. It is easy to see that 
$\lambda $ is unital, so we only check that it is multiplicative. 
We compute:
\begin{eqnarray*}
\lambda ((a\# h)(a'\# h'))&=&
\lambda ((z^1\cdot a)(z^2h_1\cdot a')\# z^3h_2h')\\
&=&T^1[(t^1_1z^1\cdot a)(t^1_2z^2h_1\cdot a')]_{(-1)}t^2z^3_1h_{(2, 1)}h'_1\\
&&\otimes T^2\cdot [(t^1_1z^1\cdot a)(t^1_2z^2h_1\cdot a')]_{(0)}
\# T^3t^3z^3_2h_{(2, 2)}h'_2\\
&{{\rm (\ref{q3})}\atop =}&T^1[(t^1y^1\cdot a)
(t^2y^2_1x^1h_1\cdot a')]_{(-1)}
t^3y^2_2x^2h_{(2, 1)}h'_1\\
&&\otimes T^2\cdot [(t^1y^1\cdot a)(t^2y^2_1x^1h_1\cdot a')]_{(0)}
\# T^3y^3x^3h_{(2, 2)}h'_2\\
&{{\rm (\ref{multi}, \;\ref{q1})}\atop =}&T^1X^1(z^1y^1\cdot a)_{(-1)}z^2
(y^2_1h_{(1, 1)}x^1\cdot a')_{(-1)}y^2_2h_{(1, 2)}x^2h'_1\\
&&\otimes  
[T^2_1X^2\cdot (z^1y^1\cdot a)_{(0)}][T^2_2X^3z^3\cdot 
(y^2_1h_{(1, 1)}x^1\cdot a')_{(0)}]\# T^3y^3h_2x^3h'_2\\
&{{\rm (\ref{yd3})}\atop =}&T^1X^1(z^1y^1\cdot a)_{(-1)}z^2
y^2_1h_{(1, 1)}(x^1\cdot a')_{(-1)}x^2h'_1\\
&&\otimes  
[T^2_1X^2\cdot (z^1y^1\cdot a)_{(0)}][T^2_2X^3z^3y^2_2h_{(1, 2)}\cdot  
(x^1\cdot a')_{(0)}]\# T^3y^3h_2x^3h'_2\\
&{{\rm (\ref{q3})}\atop =}&T^1X^1(y^1_1z^1\cdot a)_{(-1)}y^1_2z^2Z^1
h_{(1, 1)}(x^1\cdot a')_{(-1)}x^2h'_1\\
&&\otimes  
[T^2_1X^2\cdot (y^1_1z^1\cdot a)_{(0)}][T^2_2X^3y^2z^3_1Z^2h_{(1, 2)}\cdot  
(x^1\cdot a')_{(0)}]\\
&&\# T^3y^3z^3_2Z^3h_2x^3h'_2\\
&{{\rm (\ref{yd3})}\atop =}&T^1X^1y^1_1(z^1\cdot a)_{(-1)}z^2Z^1
h_{(1, 1)}(x^1\cdot a')_{(-1)}x^2h'_1\\
&&\otimes  
[T^2_1X^2y^1_2\cdot (z^1\cdot a)_{(0)}][T^2_2X^3y^2z^3_1Z^2h_{(1, 2)}\cdot  
(x^1\cdot a')_{(0)}]\\
&&\# T^3y^3z^3_2Z^3h_2x^3h'_2\\
&{{\rm (\ref{q1}, \;\ref{q3})}\atop =}&T^1(z^1\cdot a)_{(-1)}z^2
h_1Z^1(x^1\cdot a')_{(-1)}x^2h'_1\otimes   
[y^1T^2\cdot (z^1\cdot a)_{(0)}]\\
&&[y^2T^3_1z^3_1h_{(2, 1)}Z^2\cdot   
(x^1\cdot a')_{(0)}]\# y^3T^3_2z^3_2h_{(2, 2)}Z^3x^3h'_2\\
&=&[T^1(z^1\cdot a)_{(-1)}z^2h_1\otimes (T^2\cdot (z^1\cdot a)_{(0)}
\# T^3z^3h_2)]\\
&&[Z^1(x^1\cdot a')_{(-1)}x^2h'_1\otimes (Z^2\cdot (x^1\cdot a')_{(0)}
\# Z^3x^3h'_2)]\\
&=&\lambda (a\# h)\lambda (a'\# h'), \;\;\;\;\;q.e.d.
\end{eqnarray*}
Now we check (\ref{lca1}). We compute:\\[2mm]
${\;\;\;\;\;}$
$(id \otimes \lambda )(\lambda (a\# h))\Phi _{\lambda }$
\begin{eqnarray*}
&=&[T^1(t^1\cdot a)_{(-1)}t^2h_1\otimes X^1(x^1T^2\cdot (t^1\cdot a)_{(0)})
_{(-1)}x^2T^3_1t^3_1h_{(2, 1)}\\
&&\otimes X^2\cdot (x^1T^2\cdot (t^1\cdot a)_{(0)})_{(0)}\# 
X^3x^3T^3_2t^3_2h_{(2, 2)}]\Phi _{\lambda }\\
&{{\rm (\ref{q3})}\atop =}&[Y^1Z^1x^1_1(t^1\cdot a)_{(-1)}t^2h_1\otimes 
X^1(Y^2_1Z^2x^1_2\cdot (t^1\cdot a)_{(0)})
_{(-1)}Y^2_2Z^3x^2t^3_1h_{(2, 1)}\\
&&\otimes X^2\cdot (Y^2_1Z^2x^1_2\cdot (t^1\cdot a)_{(0)})_{(0)}\# 
X^3Y^3x^3t^3_2h_{(2, 2)}]\Phi _{\lambda }\\
&{{\rm (\ref{yd3})}\atop =}&Y^1Z^1x^1_1(t^1\cdot a)_{(-1)}t^2h_1T^1\otimes   
X^1Y^2_1(Z^2x^1_2\cdot (t^1\cdot a)_{(0)})
_{(-1)}Z^3x^2t^3_1h_{(2, 1)}T^2\\
&&\otimes X^2Y^2_2\cdot (Z^2x^1_2\cdot (t^1\cdot a)_{(0)})_{(0)}\# 
X^3Y^3x^3t^3_2h_{(2, 2)}T^3\\
&{{\rm (\ref{q1})}\atop =}&Y^1Z^1x^1_1(t^1\cdot a)_{(-1)}t^2T^1h_{(1, 1)}
\otimes X^1Y^2_1(Z^2x^1_2\cdot (t^1\cdot a)_{(0)})
_{(-1)}Z^3x^2t^3_1T^2h_{(1, 2)}\\
&&\otimes X^2Y^2_2\cdot (Z^2x^1_2\cdot (t^1\cdot a)_{(0)})_{(0)}\# 
X^3Y^3x^3t^3_2T^3h_2\\
&{{\rm (\ref{q3})}\atop =}&Y^1Z^1x^1_1(T^1_1z^1t^1\cdot a)_{(-1)}T^1_2z^2t^2_1 
h_{(1, 1)}\otimes   
X^1Y^2_1(Z^2x^1_2\cdot (T^1_1z^1t^1\cdot a)_{(0)})_{(-1)}\\
&&Z^3x^2T^2z^3t^2_2h_{(1, 2)}
\otimes X^2Y^2_2\cdot (Z^2x^1_2\cdot (T^1_1z^1t^1\cdot a)_{(0)})_{(0)}\#  
X^3Y^3x^3T^3t^3h_2\\
&{{\rm (\ref{yd3})}\atop =}&Y^1Z^1(z^1t^1\cdot a)_{(-1)}z^2t^2_1
h_{(1, 1)}\otimes   
X^1Y^2_1(Z^2\cdot (z^1t^1\cdot a)_{(0)})_{(-1)}
Z^3z^3t^2_2h_{(1, 2)}\\
&&\otimes X^2Y^2_2\cdot (Z^2\cdot (z^1t^1\cdot a)_{(0)})_{(0)}\#   
X^3Y^3t^3h_2\\
&{{\rm (\ref{yd1})}\atop =}&Y^1Z^1(t^1\cdot a)_{(-1)_1}t^2_1h_{(1, 1)}
\otimes X^1Y^2_1Z^2(t^1\cdot a)_{(-1)_2}t^2_2h_{(1, 2)}\\
&&\otimes X^2Y^2_2Z^3\cdot (t^1\cdot a)_{(0)}\# X^3Y^3t^3h_2\\
&{{\rm (\ref{q3})}\atop =}&X^1T^1_1(t^1\cdot a)_{(-1)_1}t^2_1h_{(1, 1)}
\otimes X^2T^1_2(t^1\cdot a)_{(-1)_2}t^2_2h_{(1, 2)}
\otimes X^3_1T^2\cdot (t^1\cdot a)_{(0)}\\
&&\# X^3_2T^3t^3h_2\\
&=&\Phi _{\lambda }(\Delta \otimes id)(\lambda (a\# h)), \;\;\;\;\;q.e.d.
\end{eqnarray*}
The conditions (\ref{lca2}), (\ref{lca3}) and (\ref{lca4}) are easy to check 
and are left to the reader.
\end{proof}
\begin{definition}
Let $H$ be a quasi-bialgebra and $(\mf {B}, \lambda , \Phi _{\lambda })$,  
$(\mf {C}, \mu , \Phi _{\mu })$ two left $H$-comodule algebras. A 
morphism of comodule algebras from $\mf {B}$ to $\mf {C}$ is an 
algebra map $\nu :\mf {B}\rightarrow  \mf {C}$ such that 
$\Phi _{\mu }=(id_H\ot id_H\ot \nu )(\Phi _{\lambda })$ and  
$\mu \circ \nu =(id_H\ot \nu )\circ \lambda $.
\end{definition}
\begin{remark}
If $H$ is a quasi-bialgebra, $A$ is an algebra in $_H^H{\cal YD}$ 
and we consider  
the left $H$-comodule algebra $A\# H$ as in Proposition \ref{smlca}, 
then one can easily see that the map $j:H\rightarrow A\# H$, 
$\;j(h)=1\# h$, is a morphism of left $H$-comodule algebras. 
\end{remark}
${\;\;\;}$
The next result is a generalization of the fact, proved in \cite{bn}, 
Proposition 2.5, that $H_0$ is an algebra in $_H^H{\cal YD}$; the proof 
is similar to the one in \cite{bn} and will be omitted. 
\begin{proposition} \label{lulu}
If $H$ is a quasi-Hopf algebra, $(\mf {B}, \lambda , \Phi _{\lambda })$  
a left $H$-comodule algebra and $v:H\rightarrow \mf {B}$ a  
morphism of $H$-comodule algebras, then $\mf {B}^v$ becomes an 
algebra in $_H^H{\cal YD}$ with coaction 
\begin{eqnarray*}
&&\mf {B}^v\rightarrow H\ot \mf {B}^v, \;\;\;\;
\mf {b}\mapsto X^1Y^1_1\mf {b}_{[-1]}g^1S(q^2Y^2_2)Y^3\ot 
v(X^2Y^1_2)\mf {b}_{[0]}v(g^2S(X^3q^1Y^2_1)),
\end{eqnarray*}
where $f^{-1}=g^1\ot g^2$ is given by (\ref{g}) and $q_R=q^1\ot q^2=
Z^1\ot S^{-1}(\alpha Z^3)Z^2$. Moreover, the map $v:H_0\rightarrow 
\mf {B}^v$ is a morphism of algebras in $_H^H{\cal YD}$.
\end{proposition}
${\;\;\;}$
As a consequence of these results, we obtain:
\begin{corollary} \label{lala}
Let $H$ be a quasi-Hopf algebra and $A$ an algebra in $_H^H{\cal YD}$. 
Then $(A\# H)^j$ becomes an algebra in $_H^H{\cal YD}$, with coaction 
$\lambda _{(A\# H)^j}:(A\# H)^j\rightarrow H\ot (A\# H)^j$, 
\begin{eqnarray*}
\lambda _{(A\# H)^j}(a\# h)&=&
X^1Y^1_1T^1(t^1\cdot a)_{(-1)}t^2h_1g^1S(q^2Y^2_2)Y^3\ot 
[X^2_1Y^1_{(2, 1)}T^2\cdot (t^1\cdot a)_{(0)}\\
&&\# X^2_2Y^1_{(2, 2)}T^3t^3h_2g^2S(X^3q^1Y^2_1)], 
\end{eqnarray*}
and the map $j:H_0\rightarrow (A\# H)^j$ is a morphism of algebras in 
$_H^H{\cal YD}$.
\end{corollary}
\begin{lemma}
Let $H$ be a quasi-Hopf algebra and $A$ an algebra in $_H^H{\cal YD}$. Then 
$i_0$ is a morphism of algebras in $_H^H{\cal YD}$ from $A$ to 
$(A\# H)^j$.
\end{lemma}
\begin{proof}
We already know that $i_0$ is a morphism of left $H$-module algebras from 
$A$ to $(A\# H)^j$, so the only thing left to prove is that 
$\lambda _{(A\# H)^j}\circ i_0=(id_H\ot i_0)\circ \lambda _A$. We first record 
the following relations, whose proofs are easy and left to the reader, 
which hold for the elements $p_R=p^1\ot p^2$ and $q_R=q^1\ot q^2$ given by 
(\ref{psiq}): 
\begin{eqnarray}
&&T^1p^1_1\ot T^2p^1_2\ot T^3p^2=y^1\ot y^2_1p^1\ot y^2_2p^2S(y^3),
\label{relatiep} \\
&&Y^1Z^1\ot q^1Y^2_1Z^2\ot S(q^2Y^2_2Z^3)Y^3=q^1_1\ot q^1_2\ot S(q^2). 
\label{relatieq}
\end{eqnarray}
Now, by denoting $p_R=P^1\ot P^2$ another copy of $p_R$, we compute: 
\begin{eqnarray*}
(\lambda _{(A\# H)^j}\circ i_0)(a)&=&\lambda _{(A\# H)^j}(p^1\cdot a
\# p^2)\\
&=&X^1Y^1_1T^1(t^1p^1\cdot a)_{(-1)}t^2p^2_1g^1S(q^2Y^2_2)Y^3\\
&&\ot [X^2_1Y^1_{(2, 1)}T^2\cdot (t^1p^1\cdot a)_{(0)}
\# X^2_2Y^1_{(2, 2)}T^3t^3p^2_2g^2S(X^3q^1Y^2_1)]\\
&{{\rm (\ref{relmare})}\atop =}&X^1Y^1_1T^1(Z^1_{(1, 1)}p^1_1P^1
\cdot a)_{(-1)}Z^1_{(1, 2)}p^1_2P^2S(q^2Y^2_2Z^3)Y^3\\
&&\ot [X^2_1Y^1_{(2, 1)}T^2\cdot (Z^1_{(1, 1)}p^1_1P^1\cdot a)_{(0)}
\# X^2_2Y^1_{(2, 2)}T^3Z^1_2p^2S(X^3q^1Y^2_1Z^2)]\\
&{{\rm (\ref{yd3}, \;\ref{q1})}\atop =}&X^1T^1Y^1_{(1, 1)}Z^1_{(1, 1)}p^1_1
(P^1\cdot a)_{(-1)}P^2S(q^2Y^2_2Z^3)Y^3\\
&&\ot [X^2_1T^2Y^1_{(1, 2)}Z^1_{(1, 2)}p^1_2\cdot (P^1\cdot a)_{(0)}
\# X^2_2T^3Y^1_2Z^1_2p^2S(X^3q^1Y^2_1Z^2)]\\
&{{\rm (\ref{relatieq})}\atop =}&X^1T^1q^1_{(1, 1, 1)}p^1_1
(P^1\cdot a)_{(-1)}P^2S(q^2)\\
&&\ot [X^2_1T^2q^1_{(1, 1, 2)}p^1_2\cdot (P^1\cdot a)_{(0)}
\# X^2_2T^3q^1_{(1, 2)}p^2S(q^1_2)S(X^3)]\\
&{{\rm (\ref{tpqr1})}\atop =}&X^1T^1p^1_1q^1_1
(P^1\cdot a)_{(-1)}P^2S(q^2)\\
&&\ot [X^2_1T^2p^1_2q^1_2\cdot (P^1\cdot a)_{(0)}
\# X^2_2T^3p^2S(X^3)]\\
&{{\rm (\ref{yd3}, \;\ref{relatiep})}\atop =}&X^1y^1(q^1_1P^1\cdot a)_{(-1)}
q^1_2P^2S(q^2)\\
&&\ot [X^2_1y^2_1p^1\cdot (q^1_1P^1\cdot a)_{(0)}
\# X^2_2y^2_2p^2S(X^3y^3)]\\
&{{\rm (\ref{tpqr2})}\atop =}&a_{(-1)}\ot (p^1\cdot a_{(0)}\# p^2)\\
&=&((id_H\ot i_0)\circ \lambda _A)(a),
\end{eqnarray*}
finishing the proof.
\end{proof}
\begin{theorem} \label{izomyd}
Let $H$ be a quasi-Hopf algebra and $A$ an algebra in $_H^H{\cal YD}$. 
Then the map 
$\Pi :A\underline{\ot }H_0\rightarrow (A\# H)^j$,  
$\;\Pi (a\ot h)=x^1\cdot a\# x^2hS(x^3)$, 
is an isomorphism of algebras in $_H^H{\cal YD}$. 
\end{theorem}  
\begin{proof}
We proved that $j:H_0\rightarrow (A\# H)^j$ and 
$i_0:A\rightarrow (A\# H)^j$ are morphisms of algebras in $_H^H{\cal YD}$, 
and together with the commutation relation (\ref{comutare}) this allows 
to apply the universal property of the twisted tensor product 
$A\underline{\ot }H_0=A\ot _RH_0$ in the category $_H^H{\cal YD}$, 
obtaining thus a morphism of algebras in $_H^H{\cal YD}$ between 
$A\underline{\ot }H_0$ and $(A\# H)^j$, $\omega $ say, which has to be given 
by $\omega (a\ot h)=i_0(a)*j(h)$, that is $\omega $ coincides 
with the map $\Pi $, which finishes the proof.
\end{proof}
\section{An algebra isomorphism $B^v\# H\simeq B\otimes H$ and applications} 
\setcounter{equation}{0}
\begin{proposition}\label{pregat}
Let $H$ be a quasi-Hopf algebra, $B$ an associative algebra and 
$v:H\rightarrow B$ an  
algebra map. Denote by $\eta $ the algebra map $\eta :H\rightarrow 
B\otimes H$, $\eta (h)=v(h_1)\otimes h_2$. Define the map
\begin{eqnarray} 
&&u:B\rightarrow B\otimes H, \;\;\;\; u(b)=v(x^1)bv(S(x^3_2X^3)f^1)
\otimes x^2X^1\beta S(x^3_1X^2)f^2,  \label{uala}
\end{eqnarray}
where $f=f^1\otimes f^2$ is the Drinfeld twist given by (\ref{f}). Then 
$u$ is a morphism of left $H$-module algebras from $B^v$ to 
$(B\otimes H)^{\eta}$. 
\end{proposition} 
\begin{proof}
The fact that $u(v(\beta ))=\eta (\beta )$ follows immediately from 
(\ref{moref}). We check now that $u$ is a morphism of left $H$-modules:
\begin{eqnarray*}
h\triangleright _{\eta }u(b)&=&\eta (h_1)u(b)\eta (S(h_2))\\
&=&v(h_{(1, 1)})v(x^1)bv(S(x^3_2X^3)f^1)v(S(h_2)_1)\otimes 
h_{(1, 2)}x^2X^1\beta S(x^3_1X^2)f^2S(h_2)_2\\
&{{\rm (\ref{ca})}\atop =}&v(h_{(1, 1)}x^1)b
v(S(h_{(2, 2)}x^3_2X^3)f^1)\otimes   
h_{(1, 2)}x^2X^1\beta S(h_{(2, 1)}x^3_1X^2)f^2\\
&{{\rm (\ref{q1})}\atop =}&v(x^1h_1)bv(S(x^3_2h_{(2, 2, 2)}X^3)f^1)\otimes    
x^2h_{(2, 1)}X^1\beta S(x^3_1h_{(2, 2, 1)}X^2)f^2\\
&{{\rm (\ref{q1},\; \ref{q5})}\atop =}&v(x^1h_1)bv(S(x^3_2X^3h_2)f^1)
\otimes  x^2X^1\beta S(x^3_1X^2)f^2\\
&=&u(h\triangleright _vb), \;\;\;\;q.e.d.
\end{eqnarray*}
Now we check that $u$ is multiplicative (we denote by $F=F^1\otimes F^2$ 
another copy of $f$):
\begin{eqnarray*}
u(b)*u(b')&=&\eta (Y^1)u(b)\eta (S(y^1Y^2)\alpha y^2Y^3_1)
u(b')\eta (S(y^3Y^3_2))\\
&=&v(Y^1_1x^1)bv(S(x^3_2X^3)f^1S(y^1Y^2)_1\alpha _1y^2_1Y^3_{(1, 1)}z^1)b'\\
&&v(S(z^3_2Z^3)F^1S(y^3Y^3_2)_1)
\otimes Y^1_2x^2X^1\beta S(x^3_1X^2)f^2S(y^1Y^2)_2\alpha _2\\
&&y^2_2Y^3_{(1, 2)}z^2Z^1\beta S(z^3_1Z^2)F^2S(y^3Y^3_2)_2\\
&{{\rm (\ref{ca},\; \ref{moref})}\atop =}&v(Y^1_1x^1)b
v(S(T^2t^1_2y^1_2Y^2_2x^3_2X^3)\alpha T^3t^2y^2_1Y^3_{(1, 1)}z^1)
b'v(S(y^3_2Y^3_{(2, 2)}z^3_2Z^3)f^1)\\
&&\otimes Y^1_2x^2X^1\beta S(T^1t^1_1y^1_1Y^2_1x^3_1X^2)\alpha t^3y^2_2
Y^3_{(1, 2)}z^2Z^1\beta S(y^3_1Y^3_{(2, 1)}z^3_1Z^2)f^2\\
&{{\rm (\ref{q1}, \;\ref{q5})}\atop =}&v(Y^1_1x^1)b
v(S(T^2t^1_2y^1_2Y^2_2x^3_2X^3)\alpha T^3t^2y^2_1z^1Y^3_1)
b'v(S(y^3_2z^3_2Z^3Y^3_2)f^1)\\
&&\otimes Y^1_2x^2X^1\beta S(T^1t^1_1y^1_1Y^2_1x^3_1X^2)\alpha t^3y^2_2
z^2Z^1\beta S(y^3_1z^3_1Z^2)f^2\\
&{{\rm (\ref{q3}, \;\ref{q1}, \;\ref{q5})}\atop =}&v(Y^1_1x^1)b
v(S(T^2t^1_2Y^2_2x^3_2X^3)\alpha T^3t^2Y^3_1)
b'v(S(y^3_2Z^3t^3Y^3_2)f^1)\\
&&\otimes Y^1_2x^2X^1\beta S(y^1T^1t^1_1Y^2_1x^3_1X^2)\alpha y^2
Z^1\beta S(y^3_1Z^2)f^2\\
&{{\rm (\ref{q3}, \;\ref{q1})}\atop =}&v(x^1Y^1)b
v(S(T^2x^3_{(1, 1, 2)}t^1_2W^2_2Y^2_{(2, 2)}X^3)\alpha T^3x^3_{(1, 2)}
t^2W^3_1Y^3_1)b'\\
&&v(S(y^3_2Z^3x^3_2t^3W^3_2Y^3_2)f^1)
\otimes x^2W^1Y^2_1X^1\beta \\
&&S(y^1T^1x^3_{(1, 1, 1)}t^1_1W^2_1
Y^2_{(2, 1)}X^2)\alpha y^2Z^1\beta S(y^3_1Z^2)f^2\\
&{{\rm (\ref{q1}, \;\ref{q5})}\atop =}&v(x^1Y^1)b
v(S(T^2t^1_2W^2_2Y^2_{(2, 2)}X^3)\alpha T^3
t^2W^3_1Y^3_1)b'v(S(y^3_2Z^3x^3_2t^3W^3_2Y^3_2)f^1)\\
&&\otimes x^2W^1Y^2_1X^1\beta S(y^1x^3_1T^1t^1_1W^2_1
Y^2_{(2, 1)}X^2)\alpha y^2Z^1\beta S(y^3_1Z^2)f^2\\
&{{\rm (\ref{q1}, \;\ref{q5})}\atop =}&v(x^1Y^1)b
v(S(T^2t^1_2W^2_2X^3Y^2)\alpha T^3
t^2W^3_1Y^3_1)b'v(S(y^3_2Z^3x^3_2t^3W^3_2Y^3_2)f^1)\\
&&\otimes x^2W^1X^1\beta S(y^1x^3_1T^1t^1_1W^2_1
X^2)\alpha y^2Z^1\beta S(y^3_1Z^2)f^2\\
&{{\rm (\ref{q3}, \;\ref{q5})}\atop =}&v(x^1Y^1)b
v(S(T^2t^1_2z^2X^3_1Y^2)\alpha T^3
t^2z^3_1X^3_{(2, 1)}Y^3_1)b'\\
&&v(S(y^3_2Z^3x^3_2t^3z^3_2X^3_{(2, 2)}Y^3_2)f^1)\\
&&\otimes x^2X^1\beta S(y^1x^3_1T^1t^1_1z^1
X^2)\alpha y^2Z^1\beta S(y^3_1Z^2)f^2\\
&{{\rm (\ref{q3}, \;\ref{q5})}\atop =}&v(x^1Y^1)b
v(S(T^2t^1_2z^2X^3_1Y^2)\alpha T^3
t^2z^3_1X^3_{(2, 1)}Y^3_1)b'\\
&&v(S(x^3_2t^3z^3_2X^3_{(2, 2)}Y^3_2)f^1)\otimes 
x^2X^1\beta S(x^3_1T^1t^1_1z^1X^2)S(y^1)\alpha y^2\beta S(y^3)f^2\\
&{{\rm (\ref{q6}, \;\ref{q3})}\atop =}&v(x^1Y^1)b
v(S(y^2_1t^1X^3_1Y^2)\alpha y^2_2t^2 
X^3_{(2, 1)}Y^3_1)b'v(S(x^3_2y^3t^3X^3_{(2, 2)}Y^3_2)f^1)\\
&&\otimes x^2X^1\beta S(x^3_1y^1X^2)f^2\\
&{{\rm (\ref{q1}, \;\ref{q5})}\atop =}&v(x^1Y^1)b
v(S(t^1Y^2)\alpha t^2Y^3_1)b'v(S(x^3_2X^3t^3Y^3_2)f^1)
\otimes x^2X^1\beta S(x^3_1X^2)f^2\\
&=&u(v(Y^1)b
v(S(t^1Y^2)\alpha t^2Y^3_1)b'v(S(t^3Y^3_2)))\\
&=&u(b*b'), 
\end{eqnarray*}
finishing the proof.
\end{proof}
\begin{corollary}
The map $\psi :B^v\# H\rightarrow B\otimes H$, given by 
\begin{eqnarray}
&&\psi (b\# h)=v(X^1x^1_1)bv(S(X^2x^1_2)\alpha X^3x^2h_1)\otimes x^3h_2, 
\label{functia1}
\end{eqnarray}
is an algebra map.
\end{corollary}
\begin{proof}
The universal property of the smash product, applied to the maps 
$\eta $ and $u$ from Proposition \ref{pregat}, provides an algebra map 
$B^v\# H\rightarrow B\otimes H$, which turns out to be exactly the map 
$\psi $ given by (\ref{functia1}), we leave the details to the reader. 
\end{proof}
\begin{proposition}
If $H$ is a quasi-Hopf algebra, $B$ an algebra, $v:H\rightarrow B$ an 
algebra map, then:\\[2mm]
(i) The map $\theta :B\rightarrow B^v\# H$, $\theta (b)=v(z^1)b
v(Z^1\beta S(z^2Z^2))\# z^3Z^3$, is an algebra map.\\
(ii) The map $\mu :H\rightarrow B^v\# H$, $\mu (h)=v(z^1Z^1\beta 
S(z^2h_1Z^2))\# z^3h_2Z^3$, is an algebra map. \\
(iii) For all $h\in H$ and $b\in B$, the following relation holds:
\begin{eqnarray*}
&&\theta (b)\mu (h)=\mu (h) \theta (b)=v(z^1)bv(Z^1\beta S(z^2h_1Z^2))
\# z^3h_2Z^3.
\end{eqnarray*}
(iv) Consequently, the map 
\begin{eqnarray}
&&\xi :B\otimes H\rightarrow B^v\# H, \;\;\;\; 
\xi (b\otimes h)=v(z^1)bv(Z^1\beta S(z^2h_1Z^2))\# z^3h_2Z^3, 
\label{functia2}
\end{eqnarray}
is an algebra map.
\end{proposition} 
\begin{proof}
We only prove (i) and leave the rest to the reader. Obviously $\theta (1)=
v(\beta )\# 1$, so we only have to check that $\theta $ is multiplicative. 
We compute:
\begin{eqnarray*}
\theta (b)\theta (b')&=&[x^1\triangleright _v(v(z^1)bv(Z^1\beta S(z^2Z^2)))]
*[x^2z^3_1Z^3_1\triangleright _v(v(t^1)b'v(T^1\beta S(t^2T^2)))]\\
&&\# x^3z^2_2Z^3_2t^3T^3\\
&=&v(X^1x^1_1z^1)bv(Z^1\beta S(y^1X^2x^1_2z^2Z^2)\alpha y^2X^3_1
x^2_1z^3_{(1, 1)}Z^3_{(1, 1)}t^1)\\
&&b'v(T^1\beta 
S(y^3X^3_2x^2_2z^3_{(1, 2)}Z^3_{(1, 2)}t^2T^2))
\# x^3z^2_2Z^3_2t^3T^3\\
&{{\rm (\ref{q3}, \;\ref{q5})}\atop =}&v(X^1w^1_1z^1y^1)bv(Z^1\beta  
S(X^2w^1_2z^2y^2_1x^1Z^2)\alpha X^3
w^2z^3_1y^2_{(2, 1)}x^2_1Z^3_{(1, 1)}t^1)\\
&&b'v(T^1\beta 
S(w^3z^3_2y^2_{(2, 2)}x^2_2Z^3_{(1, 2)}t^2T^2))
\# y^3x^3Z^3_2t^3T^3\\
&{{\rm (\ref{q3}, \;\ref{q5})}\atop =}&v(y^1)bv(Z^1\beta  
S(z^1y^2_1x^1Z^2)\alpha z^2y^2_{(2, 1)}x^2_1Z^3_{(1, 1)}t^1)\\
&&b'v(T^1\beta 
S(z^3y^2_{(2, 2)}x^2_2Z^3_{(1, 2)}t^2T^2))
\# y^3x^3Z^3_2t^3T^3\\
&{{\rm (\ref{q1}, \;\ref{q5})}\atop =}&v(y^1)bv(Z^1\beta  
S(z^1x^1Z^2)\alpha z^2x^2_1Z^3_{(1, 1)}t^1)
b'v(T^1\beta 
S(y^2z^3x^2_2Z^3_{(1, 2)}t^2T^2))\\
&&\# y^3x^3Z^3_2t^3T^3\\
&{{\rm (\ref{q3}, \;\ref{q5})}\atop =}&v(y^1)bv(Y^1Z^1\beta   
S(z^1Y^2_1Z^2)\alpha z^2Y^2_{(2, 1)}Z^3_1t^1)
b'v(T^1\beta 
S(y^2z^3Y^2_{(2, 2)}Z^3_2t^2T^2))\\
&&\# y^3Y^3t^3T^3\\
&{{\rm (\ref{q1}, \;\ref{q5})}\atop =}&v(y^1)bv(Y^1Z^1\beta   
S(z^1Z^2)\alpha z^2Z^3_1t^1)
b'v(T^1\beta S(y^2Y^2z^3Z^3_2t^2T^2))\# y^3Y^3t^3T^3\\
&{{\rm (\ref{q3}, \;\ref{q5})}\atop =}&v(y^1)bv(Y^1Z^1\beta   
S(Z^2)\alpha Z^3t^1)
b'v(T^1\beta S(y^2Y^2t^2T^2))\# y^3Y^3t^3T^3\\
&{{\rm (\ref{q6})}\atop =}&v(y^1)b
b'v(T^1\beta S(y^2T^2))\# y^3T^3\\
&=&\theta (bb'),
\end{eqnarray*}
finishing the proof.
\end{proof}
\begin{theorem} \label{teo}
The maps $\psi $ and $\xi $ given by (\ref{functia1}) and respectively 
(\ref{functia2}) are inverse to each other, providing thus an algebra 
isomorphism $B^v\# H\simeq B\otimes H$.
\end{theorem}
\begin{proof}
We compute:
\begin{eqnarray*}
\xi (\psi (b\# h))&=&v(z^1X^1x^1_1)bv(S(X^2x^1_2)\alpha X^3x^2h_1Z^1\beta 
S(z^2x^3_1h_{(2, 1)}Z^2))\# z^3x^3_2h_{(2, 2)}Z^3\\
&{{\rm (\ref{q1}, \;\ref{q5})}\atop =}&v(z^1X^1x^1_1)bv(S(X^2x^1_2)
\alpha X^3x^2Z^1
\beta S(z^2x^3_1Z^2))\# z^3x^3_2Z^3h\\
&{{\rm (\ref{q3}, \;\ref{q5})}\atop =}&v(z^1X^1Y^1_{(1, 1)}t^1_1)b
v(S(X^2Y^1_{(1, 2)}t^1_2)\alpha X^3Y^1_2t^2
\beta S(z^2Y^2t^3))\# z^3Y^3h\\
&{{\rm (\ref{q1}, \;\ref{q5})}\atop =}&v(X^1t^1_1)b
v(S(X^2t^1_2)\alpha X^3t^2
\beta S(t^3))\# h\\
&{{\rm (\ref{q3}, \;\ref{q5})}\atop =}&bv(S(t^1)\alpha t^2\beta S(t^3))\# h\\
&{{\rm (\ref{q6})}\atop =}&b\# h,
\end{eqnarray*}
\begin{eqnarray*}
\psi (\xi (b\otimes h))&=&v(X^1x^1_1z^1)bv(Z^1\beta S(X^2x^1_2z^2h_1Z^2)
\alpha X^3x^2z^3_1h_{(2, 1)}Z^3_1)\otimes x^3z^3_2h_{(2, 2)}Z^3_2\\
&{{\rm (\ref{q3}, \;\ref{q5})}\atop =}&bv(Z^1\beta S(t^1h_1Z^2)
\alpha t^2h_{(2, 1)}Z^3_1)\otimes t^3h_{(2, 2)}Z^3_2\\
&{{\rm (\ref{q1}, \;\ref{q5})}\atop =}&bv(Z^1\beta S(t^1Z^2)
\alpha t^2Z^3_1)\otimes ht^3Z^3_2\\
&{{\rm (\ref{q3}, \;\ref{q5})}\atop =}&bv(Z^1\beta S(Z^2)
\alpha Z^3)\otimes h\\
&{{\rm (\ref{q6})}\atop =}&b\otimes h,
\end{eqnarray*}
finishing the proof.
\end{proof}
${\;\;\;}$We present now some applications of Theorem \ref{teo}. We first  
record two results of independent interest, whose proofs are  
straightforward and will be omitted.  
\begin{lemma} \label{unu}
If $H$ is a quasi-Hopf algebra, then $(H\otimes H, \lambda , 
\Phi _{\lambda })$ is a left $H$-comodule algebra, with structures 
\begin{eqnarray*}
&&\lambda :H\otimes H\rightarrow H\otimes (H\otimes H), \;\;\;
\lambda (h'\otimes h)=\Phi \cdot ((\Delta \otimes id_H)(h'\otimes h))\cdot 
\Phi ^{-1}, \\
&&\Phi _{\lambda }=(id_H\otimes id_H\otimes \Delta )(\Phi )\in H\otimes H
\otimes (H\otimes H).
\end{eqnarray*}
Moreover, $\Delta :H\rightarrow H\otimes H$ is a morphism of left 
$H$-comodule algebras.
\end{lemma}
\begin{lemma} \label{doi}
If $H$ is a quasi-Hopf algebra, $\nu :\mf {B}\rightarrow  \mf {C}$ is a 
morphism of left $H$-comodule algebras and $v:H\rightarrow \mf {B}$, 
$w:H\rightarrow \mf {C}$ are morphisms of left $H$-comodule algebras 
such that $\nu \circ v=w$, then $\nu $ is a morphism of algebras in 
$_H^H{\cal YD}$ between $\mf {B}^v$ and $\mf {C}^w$, where 
$\mf {B}^v$ and $\mf {C}^w$ become algebras in $_H^H{\cal YD}$ as in 
Proposition \ref{lulu}.  
\end{lemma}
${\;\;\;}$We consider now the particular case $B=H$ and $v=id_H$ in 
Theorem \ref{teo}, so we have an algebra isomorphism 
\begin{eqnarray}
&&\Psi :H_0\# H\rightarrow H\otimes H, \;\;\;
\Psi (h'\otimes h)=X^1x^1_1h'S(X^2x^1_2)\alpha X^3x^2h_1\otimes x^3h_2. 
\label{psipartic}
\end{eqnarray}
By a slightly longer but also straightforward computation, one can 
prove the following result:
\begin{lemma}
If we consider $H_0\# H$ a left $H$-comodule algebra as in Proposition 
\ref{smlca} and $H\otimes H$ a left $H$-comodule algebra as in 
Lemma \ref{unu}, then the map $\Psi $ given by (\ref{psipartic}) is an 
isomorphism of left $H$-comodule algebras. 
\end{lemma} 
${\;\;\;}$By using (\ref{q5}) and (\ref{q6}), one can see that  
$\Psi \circ j=\Delta $, where $j:H\rightarrow H_0\# H$, $j(h)=1_{H_0}\# h$. 
Thus, as a consequence of the above results and of  
Theorem \ref{izomyd}, we finally obtain:
\begin{proposition}
If $H$ is a quasi-Hopf algebra, then $H_0\underline{\otimes }H_0\simeq 
(H\otimes H)^{\Delta }$ as algebras in $_H^H{\cal YD}$.
\end{proposition}
${\;\;\;}$As a second application of Theorem \ref{teo}, we will obtain 
a duality theorem for quasi-Hopf algebras. We recall from \cite{hn1} 
the construction of the two-sided crossed product 
$H\gsl H^*\trl H$ associated to a finite dimensional quasi-Hopf algebra $H$, 
which is an associative algebra structure on $H\otimes H^*\otimes H$ 
with multiplication 
\begin{eqnarray*}
&&(h\otimes \varphi \otimes l)(h'\otimes \varphi '\otimes l')=
hh'_1x^1\otimes (y^1\rightharpoonup \varphi \leftharpoonup h'_2x^2)
(y^2l_1\rightharpoonup \varphi '\leftharpoonup x^3)\otimes y^3l_2l',
\end{eqnarray*}
for $h, h', l, l'\in H$ and $\varphi , \varphi '\in H^*$, where 
$\rightharpoonup $ and $\leftharpoonup $ are the regular actions of $H$ 
on $H^*$. We also recall from \cite{bc} the construction of the 
so-called quasi-smash product, denoted by $H\overline{\# }H^*$, 
which is a left $H$-module algebra structure on $H\otimes H^*$, with 
multiplication and $H$-action given by 
\begin{eqnarray*}
(h\otimes \varphi )(h'\otimes \varphi ')&=&hh'_1x^1\otimes 
(\varphi \leftharpoonup h'_2x^2)(\varphi '\leftharpoonup x^3), \\
h\cdot (h'\otimes \varphi )&=&h'\otimes h\rightharpoonup \varphi ,
\end{eqnarray*}
for all $h, h'\in H$ and $\varphi , \varphi '\in H^*$. By \cite{bc}, 
the following identification of algebras holds:
\begin{eqnarray*}
&&H\gsl H^*\trl H \equiv (H\overline{\# }H^*)\# H.
\end{eqnarray*}
By \cite{pvo}, $H\overline{\# }H^*$ is isomorphic, as  
left $H$-module algebras, to $End (H)^v$, where $End (H)$ is regarded 
as an associative algebra in the usual way and $v:H\rightarrow End (H)$ 
is a certain algebra map. \\
${\;\;\;}$All these facts combined with Theorem \ref{teo} yield the 
desired duality theorem:
\begin{theorem}
If $H$ is a finite dimensional quasi-Hopf algebra, then the two-sided  
crossed product $H\gsl H^*\trl H$ is isomorphic to 
$End (H)\otimes H$ as associative algebras. 
\end{theorem} 
\section{Some iterated products}\label{sec5}
\setcounter{equation}{0}
${\;\;\;}$
Let $H$ be a quasi-bialgebra and $A$ an algebra in $_H^H{\cal YD}$. 
By Proposition \ref{smlca}, $(A\# H, \lambda , 
\Phi _{\lambda})$ is a left $H$-comodule algebra,  
and by the Preliminaries $(A\# H, \rho , \Phi _{\rho })$ is a right 
$H$-comodule algebra.      
\begin{proposition}
$(A\# H, \lambda , \rho , \Phi _{\lambda }, \Phi _{\rho }, 
\Phi _{\lambda , \rho })$ is an $H$-bicomodule algebra, with 
$\Phi _{\lambda , \rho }=X^1\ot (1\# X^2)\ot X^3\in H\ot (A\# H)\ot H$. 
\end{proposition} 
\begin{proof}
We first check the relation (\ref{bca1}). We compute: \\
${\;\;\;}$$ \Phi _{\lambda , \rho }(\lambda \ot id)(\rho (a\# h))$
\begin{eqnarray*}
&=&\Phi _{\lambda , \rho }(\lambda \ot id)((x^1\cdot a\# x^2h_1)\ot x^3h_2)\\
&=&\Phi _{\lambda , \rho }(T^1(t^1x^1\cdot a)_{(-1)}t^2x^2_1h_{(1, 1)}\ot 
(T^2\cdot (t^1x^1\cdot a)_{(0)}\# T^3t^3x^2_2h_{(1, 2)})\ot x^3h_2)\\
&{{\rm (\ref{q3})}\atop =}&X^1T^1(t^1_1x^1\cdot a)_{(-1)}t^1_2x^2Y^1
h_{(1, 1)}\ot   
(X^2_1T^2\cdot (t^1_1x^1\cdot a)_{(0)}\# X^2_2T^3t^2x^3_1Y^2h_{(1, 2)})\\
&&\ot  
X^3t^3x^3_2Y^3h_2\\
&{{\rm (\ref{q1}, \;\ref{yd3})}\atop =}&X^1T^1t^1_1(x^1\cdot a)_{(-1)}
x^2h_1Y^1\ot   
(X^2_1T^2t^1_2\cdot (x^1\cdot a)_{(0)}\# X^2_2T^3t^2x^3_1h_{(2, 1)}Y^2)\\
&&\ot  
X^3t^3x^3_2h_{(2, 2)}Y^3\\
&{{\rm (\ref{q3})}\atop =}&X^1(x^1\cdot a)_{(-1)}x^2h_1Y^1\ot (y^1X^2\cdot  
(x^1\cdot a)_{(0)}\# y^2X^3_1x^3_1h_{(2, 1)}Y^2)\\
&&\ot y^3X^3_2x^3_2h_{(2, 2)}Y^3\\
&=&(X^1(x^1\cdot a)_{(-1)}x^2h_1\ot (y^1X^2\cdot (x^1\cdot a)_{(0)}\# 
y^2X^3_1x^3_1h_{(2, 1)})\ot y^3X^3_2x^3_2h_{(2, 2)})\Phi _{\lambda ,  
\rho }\\
&=&(id \ot \rho )(X^1(x^1\cdot a)_{(-1)}x^2h_1\ot (X^2\cdot 
(x^1\cdot a)_{(0)}\# X^3x^3h_2))\Phi _{\lambda , \rho }\\
&=&(id \ot \rho )(\lambda (a\# h))\Phi _{\lambda , \rho },\;\;\; q.e.d.  
\end{eqnarray*}
Then, using the fact that $\lambda (1\# h)=h_1\ot (1\# h_2)$ and 
$\rho (1\# h)=(1\# h_1)\ot h_2$, one can see immediately that the 
conditions (\ref{bca2}) and (\ref{bca3}) reduce to the condition 
(\ref{q3}) for $\Phi $.
\end{proof}
${\;\;\;}$
Let now again $H$ be a quasi-bialgebra, $A$ an algebra in 
$_H^H{\cal YD}$ and $\cal A$ an $H$-bimodule algebra. We can consider the 
$H$-bicomodule algebra $A\# H$ as above, and then the L-R-smash product 
${\cal A}\nat (A\# H)$. On the other hand, we can first consider the 
$H$-bimodule algebra ${\cal A}\odot A$, and then the L-R-smash product 
$({\cal A}\odot A)\nat H$.  
\begin{proposition}
There is an algebra isomorphism 
$\Psi :{\cal A}\nat (A\# H)\simeq ({\cal A}\odot A)\nat H$, given by
\begin{eqnarray*} 
&&\Psi \;(\varphi \nat (a\# h))=(x^1\cdot \varphi \otimes x^2\cdot a)\nat 
x^3h, \;\;\;\forall \;\;\varphi \in {\cal A}, \;a\in A, \;h\in H.
\end{eqnarray*} 
\end{proposition}
\begin{proof}
Obviously $\Psi $ is bijective and unital, so we only have to check that 
it is multiplicative. We first write down the multiplication rule in 
${\cal A}\nat (A\# H)$:\\[2mm]
${\;\;\;\;\;}$
$[\varphi \nat (a\# h)][\varphi '\nat (a'\# h')]$
\begin{eqnarray*}
&=&(\tilde{x}^1_{\lambda }\cd \varphi \cd 
\theta ^3(a'\# h')_{<1>}\tilde{x}^2_{\rho })
(\tilde{x}^2_{\lambda }(a\# h)_{[-1]}\theta ^1\cd \varphi ' \cd 
\tilde{x}^3_{\rho })
\nat \tilde{x}^3_{\lambda }(a\# h)_{[0]}\theta ^2
(a'\# h')_{<0>}\tilde{x}^1_{\rho }\\
&=&(x^1\cdot \varphi \cdot t^3w^3h'_2y^2)(x^2Y^1(z^1\cdot a)_{(-1)}
z^2h_1t^1\cdot \varphi '\cdot y^3)\\
&&\nat (1\# x^3)(Y^2\cdot (z^1\cdot a)_{(0)}\# Y^3z^3h_2)(1\# t^2)
(w^1\cdot a'\# w^2h'_1)(1\# y^1)\\
&=&(x^1\cdot \varphi \cdot t^3w^3h'_2y^2)(x^2Y^1(z^1\cdot a)_{(-1)}
z^2h_1t^1\cdot \varphi '\cdot y^3)\\
&&\nat (v^1x^3_1Y^2\cdot (z^1\cdot a)_{(0)})(v^2x^3_{(2, 1)}Y^3_1z^3_1
h_{(2, 1)}t^2_1w^1\cdot a')\# v^3x^3_{(2, 2)}Y^3_2z^3_2h_{(2, 2)}
t^2_2w^2h'_1y^1,
\end{eqnarray*}
where for the second equality we have used the formulae, presented before, 
giving the $H$-bicomodule algebra structure of $A\# H$. Now we compute:\\[2mm]
${\;\;}$
$\Psi (\varphi \nat (a\# h))\Psi (\varphi '\nat (a'\# h'))$
\begin{eqnarray*}
&=&[(z^1\cdot \varphi \otimes z^2\cdot a)\nat z^3h][(w^1\cdot \varphi '
\otimes w^2\cdot a')\nat w^3h']\\
&=&(x^1\cdot (z^1\cdot \varphi \otimes z^2\cdot a)\cdot t^3w^3_2h'_2y^2)
(x^2z^3_1h_1t^1\cdot (w^1\cdot \varphi '\otimes w^2\cdot a')\cdot y^3)\\
&&\nat x^3z^3_2h_2t^2w^3_1h'_1y^1\\ 
&=&(x^1_1z^1\cdot \varphi \cdot t^3w^3_2h'_2y^2\otimes x^1_2z^2\cdot a)
(x^2_1z^3_{(1, 1)}h_{(1, 1)}t^1_1w^1\cdot \varphi '\cdot y^3\\
&&\otimes x^2_2z^3_{(1, 2)}h_{(1, 2)}t^1_2w^2\cdot a')
\nat x^3z^3_2h_2t^2w^3_1h'_1y^1\\ 
&=&(u^1X^1x^1_1z^1\cdot \varphi \cdot t^3w^3_2h'_2y^2)
(u^2Y^1(v^1X^2x^1_2z^2\cdot a)_{(-1)}\\
&&v^2X^3_1x^2_1z^3_{(1, 1)}h_{(1, 1)}t^1_1w^1\cdot \varphi '\cdot y^3)
\otimes (u^3_1Y^2\cdot (v^1X^2x^1_2z^2\cdot a)_{(0)})\\
&&(u^3_2Y^3v^3X^3_2x^2_2z^3_{(1, 2)}h_{(1, 2)}t^1_2w^2\cdot a')
\nat x^3z^3_2h_2t^2w^3_1h'_1y^1\\
&{{\rm (\ref{q3}, \;\ref{q1})}\atop =}&
(u^1x^1\cdot \varphi \cdot t^3w^3_2h'_2y^2)
(u^2Y^1(x^2_{(1, 1)}v^1z^1\cdot a)_{(-1)}
x^2_{(1, 2)}v^2z^2_1h_{(1, 1)}t^1_1w^1\cdot \varphi '\cdot y^3)\\
&&\otimes (u^3_1Y^2\cdot (x^2_{(1, 1)}v^1z^1\cdot a)_{(0)})
(u^3_2Y^3x^2_2v^3z^2_2h_{(1, 2)}t^1_2w^2\cdot a')
\nat x^3z^3h_2t^2w^3_1h'_1y^1\\
&{{\rm (\ref{q3}, \;\ref{yd3})}\atop =}&
(u^1x^1\cdot \varphi \cdot t^3w^3_2h'_2y^2)
(u^2Y^1x^2_{(1, 1)}v^1_1(z^1\cdot a)_{(-1)}
z^2X^1h_{(1, 1)}t^1_1w^1\cdot \varphi '\cdot y^3)\\
&&\otimes (u^3_1Y^2x^2_{(1, 2)}v^1_2\cdot (z^1\cdot a)_{(0)})
(u^3_2Y^3x^2_2v^2z^3_1X^2h_{(1, 2)}t^1_2w^2\cdot a')\\
&&\nat x^3v^3z^3_2X^3h_2t^2w^3_1h'_1y^1\\
&{{\rm (\ref{q1}, \;\ref{q3})}\atop =}&
(u^1x^1\cdot \varphi \cdot t^3w^3h'_2y^2)
(u^2Y^1x^2_{(1, 1)}v^1_1(z^1\cdot a)_{(-1)}
z^2h_1t^1\cdot \varphi '\cdot y^3)\\
&&\otimes (u^3_1Y^2x^2_{(1, 2)}v^1_2\cdot (z^1\cdot a)_{(0)})
(u^3_2Y^3x^2_2v^2z^3_1h_{(2, 1)}t^2_1w^1\cdot a')\\
&&\nat x^3v^3z^3_2h_{(2, 2)}t^2_2w^2h'_1y^1\\
&{{\rm (\ref{q1}, \;\ref{q3})}\atop =}&
(u^1_1x^1\cdot \varphi \cdot t^3w^3h'_2y^2)
(u^1_2x^2Y^1(z^1\cdot a)_{(-1)}
z^2h_1t^1\cdot \varphi '\cdot y^3)\\
&&\otimes (u^2_1v^1x^3_1Y^2\cdot (z^1\cdot a)_{(0)})
(u^2_2v^2x^3_{(2, 1)}Y^3_1z^3_1h_{(2, 1)}t^2_1w^1\cdot a')\\
&&\nat u^3v^3x^3_{(2, 2)}Y^3_2z^3_2h_{(2, 2)}t^2_2w^2h'_1y^1\\
&=&\Psi ([\varphi \nat (a\# h)][\varphi '\nat (a'\# h')]), 
\end{eqnarray*}
finishing the proof.
\end{proof}
\begin{corollary} \label{sitpartic}
If $H$ is a quasi-bialgebra, $A$ an algebra in  
$_H^H{\cal YD}$ and $C$ a left $H$-module algebra, then we have an 
algebra isomorphism  
$\Psi :C\;\gsm (A\# H)\simeq (C\diamond A)\# H$, given by  
$\Psi \;(c\;\gsm (a\# h))=(x^1\cdot c\otimes x^2\cdot a)\#   
x^3h$, for all $c\in C$, $a\in A$, $h\in H$.
\end{corollary}
${\;\;\;}$
We recall from \cite{bpvo} that, 
if $H$ is a quasi-bialgebra, $C$ a left $H$-module algebra 
and $\mb A$ an $H$-bicomodule algebra, then $C\;\gsm {\mb A}$ becomes  
a right $H$-comodule algebra, with structure 
defined for all $c\in C$ and $u\in {\mb A}$ by:    
\begin{eqnarray*}
&&\r : C\;\gsm {\mb A}\ra (C\;\gsm {\mb A})\ot H,
\mbox{${\;\;}$}
\r (c\;\gsm u)=(\theta ^1\cd c\;\gsm \theta ^2u_{<0>})
\ot \theta ^3u_{<1>},\\
&&\Phi _{\r }=(1_C\;\gsm \tX ^1_{\r })\ot \tX ^2_{\r }
\ot \tX ^3_{\r }\in (C\;\gsm {\mb A})\ot H\ot H.
\end{eqnarray*}
\begin{proposition}
The map $\Psi $ in Corollary \ref{sitpartic} is an isomorphism of right 
$H$-comodule algebras.
\end{proposition}
\begin{proof}
Note first that $C\;\gsm (A\# H)$ and $(C\diamond A)\# H$ are indeed 
right $H$-comodule algebras, since $C$ and $C\diamond A$ are left $H$-module 
algebras and $A\# H$ and $H$ are $H$-bicomodule algebras. It is easy to 
see that $\Psi $ respects the reassociators, so we only prove that 
it intertwines the coactions:\\[2mm]
${\;\;\;\;\;}$
$(\Psi \ot id_H)\circ \rho _{C\;\gsm (A\# H)}(c\;\gsm (a\# h))$
\begin{eqnarray*}
&=&(\Psi \ot id_H)((\theta ^1\cd c\;\gsm \theta ^2(a\# h)_{<0>})
\ot \theta ^3(a\# h)_{<1>})\\
&=&\Psi (z^1\cdot c\;\gsm (1\# z^2)(x^1\cdot a\# x^2h_1))\ot z^3x^3h_2\\
&=&(y^1z^1\cdot c\ot y^2z^2_1x^1\cdot a)\# y^3z^2_2x^2h_1\ot z^3x^3h_2\\
&{{\rm (\ref{q3})}\atop =}&(y^1_1x^1\cdot c\ot y^1_2x^2\cdot a)\# 
y^2x^3_1h_1\ot y^3x^3_2h_2\\
&=&\rho _{(C\diamond A)\# H}((x^1\cdot c\ot x^2\cdot a)\# x^3h)\\
&=&\rho _{(C\diamond A)\# H}\circ \Psi (c\;\gsm (a\# h)),
\end{eqnarray*}
finishing the proof.
\end{proof}
${\;\;\;}$
For our next results we will need the following slight generalization 
of Proposition \ref{izomleft}, whose proof is similar and will 
be omitted.
\begin{proposition} \label{izomleftmore}
Let $H$ be a quasi-Hopf algebra, $A$ a left $H$-module algebra, 
$\mathfrak{B}$ a left $H$-comodule algebra and $v:H\rightarrow 
\mathfrak{B}$ a morphism of left $H$-comodule algebras. Then the map 
$\tilde{v}:H\rightarrow A\;\gsm \;\mathfrak{B}$, 
$\tilde{v}(h)=1_A\;\gsm v(h)$, is an algebra map and the  
map $\Pi :A\diamond \mathfrak{B}^v\rightarrow 
(A\;\gsm \;\mathfrak{B})^{\tilde{v}}$,   
$\;\Pi (a\ot \mathfrak{b})=x^1\cdot a\;\gsm v(x^2)\mathfrak{b}v(S(x^3))$,      
is an isomorphism of left $H$-module algebras. 
\end{proposition}    
\begin{proposition}
Let $H$ be a quasi-Hopf algebra, $A$ an algebra in   
$_H^H{\cal YD}$ and $C$ a left $H$-module algebra. Then we have an 
isomorphism of algebras in $_H{\cal M}$:
\begin{eqnarray}
&&(C\diamond A)\diamond H_0\simeq C\diamond (A\underline{\otimes }H_0), 
\;\;\;(c\ot a)\ot h\mapsto X^1\cdot c\ot (X^2\cdot a\ot X^3\triangleright h). 
\label{cucu}
\end{eqnarray}
\end{proposition}
\begin{proof}
We have $(C\diamond A)\diamond H_0\simeq ((C\diamond A)\# H)^j$ by 
Proposition \ref{izomleft}, then $((C\diamond A)\# H)^j\simeq 
(C\;\gsm  (A\# H))^{\tilde{\jmath}}$ by using Corollary \ref{sitpartic}, 
then   
$(C\;\gsm  (A\# H))^{\tilde{\jmath}}\simeq C\diamond (A\# H)^j$ by   
Proposition \ref{izomleftmore} and finally  
$C\diamond (A\# H)^j\simeq C\diamond (A\diamond H_0)$ again by 
Proposition \ref{izomleft}. Moreover, one can check that the 
composition of these four isomorphisms in $_H{\cal M}$ is exactly the map 
given by (\ref{cucu}).
\end{proof}
${\;\;\;}$
We recall now some more facts from \cite{bpvo}. Let $H$ be a quasi-bialgebra, 
$A$ a left $H$-module algebra and $B$ a right $H$-module algebra. The 
two-sided smash product $A\# H\# B$ is the following associative algebra 
structure on $A\ot H\ot B$:
\begin{eqnarray*}
(a\# h\# b)(a{'}\# h{'}\# b{'})=
(x^1\cd a)(x^2h_1y^1\cd a')\# x^3h_2y^2h'_1z^1\# 
(b\cd y^3h'_2z^2)(b'\cd z^3), 
\end{eqnarray*}
for all $a, a{'}\in A$, $h, h{'}\in H$ and $b, b{'}\in B$. For the 
particular case when $A=k$, this reduces to the (right-handed) smash product 
$H\# B$, whose multiplication is 
\begin{eqnarray*}
(h\# b)(h{'}\# b{'})=hh'_1x^1\# (b\cd h'_2x^2)(b{'}\cd x^3), 
\end{eqnarray*}
for all $b, b{'}\in B$, $h, h{'}\in H$. Moreover, $H\# B$ is a left 
$H$-comodule algebra, with structures 
\begin{eqnarray*}
&&\lambda : H\# B\ra H\ot (H\# B),~~
\l (h\# b)=h_1x^1\ot (h_2x^2\# b\cd x^3),
\mbox{${\;\;}$$\forall $ $h\in H$, $b\in B$,}\\
&&\Phi _{\l }=X^1\ot X^2\ot (X^3\# 1_B)\in H\ot H\ot (H\# B),  
\end{eqnarray*}
and the map $v:H\rightarrow H\# B$, $v(h)=h\# 1$, is a morphism of 
left $H$-comodule algebras. Moreover, it was proved in \cite{bpvo} that 
with respect to this comodule algebra structure of $H\# B$ we have 
$A\# H\# B=A\;\gsm (H\# B)$. \\
${\;\;\;}$
Now, if we assume that $H$ is a quasi-Hopf algebra and if we denote by 
$j:H\rightarrow A\# H\# B$, $j(h)=1_A\# h\# 1_B$, which is an algebra map, 
as a consequence of Proposition \ref{izomleftmore} and of the above facts 
we obtain the following description of the left $H$-module algebra 
$(A\# H\# B)^j$ as a twisted tensor product:
\begin{proposition}
$(A\# H\# B)^j\simeq A\diamond (H\# B)^v$ as left $H$-module algebras. 
\end{proposition}
${\;\;\;}$We discuss now some different kinds of iterated products, 
namely in the sense of \cite{jlpvo}. We begin by recalling the framework 
in \cite{jlpvo}, extended to the case of nontrivial associators. Let 
${\cal C}$ be a monoidal category, $A$, $B$, $C$ algebras in ${\cal C}$ 
and $R_1:B\otimes A\rightarrow A\otimes B$, 
$R_2:C\otimes B\rightarrow B\otimes C$, 
$R_3:C\otimes A\rightarrow A\otimes C$ twisting maps in ${\cal C}$. Define 
the maps 
\begin{eqnarray*}
&&T_1:C\otimes (A\otimes B)\rightarrow (A\otimes B)\otimes C, \;\;\;   
T_2:(B\otimes C)\otimes A\rightarrow A\otimes (B\otimes C), \\
&&T_1=a_{A, B, C}^{-1}\circ (id_A\otimes R_2)\circ a_{A, C, B}\circ 
(R_3\otimes id_B)\circ a_{C, A, B}^{-1}, \\
&&T_2=a_{A, B, C}\circ (R_1\otimes id_C)\circ a_{B, A, C}^{-1}\circ  
(id_B\otimes R_3)\circ a_{B, C, A}.
\end{eqnarray*}
Assume that $R_1$, $R_2$, $R_3$ satisfy the braid equation in ${\cal C}$:
\begin{eqnarray}
&&(id_A\otimes R_2)\circ a_{A, C, B}\circ (R_3\otimes id_B)\circ 
a_{C, A, B}^{-1}\circ (id_C\otimes R_1)\circ a_{C, B, A} \nonumber \\
&&\hspace*{1cm}
=a_{A, B, C}\circ (R_1\otimes id_C)\circ a_{B, A, C}^{-1}\circ 
(id_B\otimes R_3)\circ a_{B, C, A}\circ (R_2\otimes id_A). \label{braid}
\end{eqnarray}
Then $T_1$ is a twisting map between $A\otimes _{R_1}B$ and $C$, $T_2$ is a 
twisting map between $A$ and $B\otimes _{R_2}C$, and moreover 
$a_{A, B, C}$ is an algebra isomorphism between 
$(A\otimes _{R_1}B)\otimes _{T_1}C$ and $A\otimes _{T_2}(B\otimes _{R_2}C)$.\\ 
${\;\;\;}$We present first a quasialgebra analogue of a result in 
\cite{jlpvo}.
\begin{proposition}
Let $H$ be a quasi-Hopf algebra, $A$, $B$ two left $H$-module algebras and 
$C$ an algebra in $^H_H{\cal YD}$. Assume that $R_1:B\otimes A\rightarrow 
A\otimes B$ is a twisting map in $_H{\cal M}$ and denote by 
$R_2:C\otimes B\rightarrow B\otimes C$ and $R_3:C\otimes A\rightarrow A
\otimes C$ the twisting maps in $_H{\cal M}$ given by 
$R_2(c\otimes b)=c_{(-1)}\cdot b\otimes c_{(0)}$ and $R_3(c\otimes a)=
c_{(-1)}\cdot a\otimes c_{(0)}$, for all $a\in A$, $b\in B$, $c\in C$. 
Then the twisting maps $R_1$, $R_2$, $R_3$ satisfy the braid equation 
(\ref{braid}) in $_H{\cal M}$ and the iterated product 
$(A\otimes _{R_1}B)\otimes _{T_1}C$ coincides with 
$(A\otimes _{R_1}B)\diamond C$.
\end{proposition}
\begin{proof}
For $a\in A$ and $b\in B$, we denote $R_1(b\otimes a)=a_{R_1}\otimes 
b_{R_1}$. We check (\ref{braid}):\\[2mm]
${\;\;\;\;}$
$(id_A\otimes R_2)\circ a_{A, C, B}\circ (R_3\otimes id_B)\circ 
a_{C, A, B}^{-1}\circ (id_C\otimes R_1)\circ a_{C, B, A}
(c\otimes b\otimes a)$
\begin{eqnarray*}
&=&Z^1(x^1Y^1\cdot c)_{(-1)}x^2\cdot (Y^3\cdot a)_{R_1}\otimes 
[Z^2\cdot (x^1Y^1\cdot c)_{(0)}]_{(-1)}Z^3x^3\cdot (Y^2\cdot b)_{R_1}\\
&&\otimes [Z^2\cdot (x^1Y^1\cdot c)_{(0)}]_{(0)}\\
&{{\rm (\ref{yd1})}\atop =}&
X^1(Y^1\cdot c)_{(-1)_1}\cdot (Y^3\cdot a)_{R_1}\otimes 
X^2(Y^1\cdot c)_{(-1)_2}\cdot (Y^2\cdot b)_{R_1}
\otimes X^3\cdot (Y^1\cdot c)_{(0)}\\
&=&a_{A, B, C}((Y^1\cdot c)_{(-1)}\cdot R_1(Y^2\cdot b\otimes Y^3\cdot a)
\otimes (Y^1\cdot c)_{(0)})\\
&=&a_{A, B, C}(R_1((Y^1\cdot c)_{(-1)_1}Y^2\cdot b\otimes 
(Y^1\cdot c)_{(-1)_2}Y^3\cdot a)\otimes (Y^1\cdot c)_{(0)})\\
&{{\rm (\ref{yd1})}\atop =}&
a_{A, B, C}(R_1(x^1X^1c_{(-1)}\cdot b\otimes x^2(X^2\cdot c_{(0)})_{(-1)}
X^3\cdot a)\otimes x^3\cdot (X^2\cdot c_{(0)})_{(0)})\\
&=&a_{A, B, C}\circ (R_1\otimes id_C)\circ a_{B, A, C}^{-1}\circ 
(id_B\otimes R_3)\circ a_{B, C, A}\circ (R_2\otimes id_A)
(c\otimes b\otimes a).
\end{eqnarray*}
A similar computation shows that in this case the map $T_1$ is given by 
$T_1(c\otimes a\otimes b)=c_{(-1)_1}\cdot a\otimes c_{(-1)_2}\cdot b
\otimes c_{(0)}$, which implies immediately that  
$(A\otimes _{R_1}B)\otimes _{T_1}C=(A\otimes _{R_1}B)\diamond C$. 
\end{proof}
${\;\;\;}$Second, we will prove that iterating the quasiassociative 
version of the Clifford process fits inside the theory developed 
in \cite{jlpvo}. Let $A$ be a left module algebra over a quasi-bialgebra $H$, 
let $\sigma $ be an $H$-linear involutive unital automorphism of $A$, 
let $q$, $s$ be two nonzero scalars and denote $C(k, q)=k[v]/(v^2=q)$ 
and $C(k, s)=k[w]/(w^2=s)$. Consider the linear maps 
$R_1:C(k, q)\otimes A\rightarrow A\otimes C(k, q)$, 
$R_2:C(k, s)\otimes C(k, q)\rightarrow C(k, q)\otimes C(k, s)$, 
$R_3:C(k, s)\otimes A\rightarrow A\otimes C(k, s)$, given by 
$R_1(1\otimes a)=a\otimes 1$, $R_1(v\otimes a)=\sigma (a)\otimes v$, 
$R_2(1\otimes 1)=1\otimes 1$, $R_2(1\otimes v)=v\otimes 1$, $R_2(w\otimes 1)=
1\otimes w$, $R_2(w\otimes v)=-v\otimes w$, $R_3(1\otimes a)=a\otimes 1$, 
$R_3(w\otimes a)=\sigma (a)\otimes w$. 
We consider $C(k, q)$ and $C(k, s)$ as left $H$-module algebras with 
trivial $H$-action, and so $R_1$, $R_2$, $R_3$ become twisting maps 
in the category $_H{\cal M}$. One can easily check that they satisfy the 
braid equation (\ref{braid}), so we have the two twisting maps $T_1$ and 
$T_2$ and moreover $(A\otimes _{R_1}C(k, q))\otimes _{T_1}C(k, s)\equiv  
A\otimes _{T_2}(C(k, q)\otimes _{R_2}C(k, s))$, because 
$a_{A, C(k, q), C(k, s)}$ is the identity. One can also check that, if we 
regard $\overline{A}=A\otimes _{R_1}C(k, q)$ as the Clifford process 
applied to the pair $(A, \sigma )$, then 
$(A\otimes _{R_1}C(k, q))\otimes _{T_1}C(k, s)$ is the Clifford process 
applied to the pair $(\overline{A}, \overline{\sigma })$, where 
$\overline{\sigma }$ is the extended automorphism given by 
$\overline{\sigma }(a\otimes 1+b\otimes v)=\sigma (a)\otimes 1-
\sigma (b)\otimes v$.\\
${\;\;\;}$Note that in the strictly associative case all these facts were 
obtained in \cite{am2} as an immediate consequence of the fact that the  
associative Clifford process is given (unlike the quasiassociative one) 
by a {\it braided} tensor product of algebras.      
\section{Universal properties}
\setcounter{equation}{0}
${\;\;\;}$We have seen in Proposition \ref{izomleft} a relation 
between the algebra $A\# H$ and the quasialgebra $A\diamond H_0$. We 
analyze now the relation between the universal property of $A\# H$ and 
the universal property of $A\diamond H_0$ regarded as a twisted tensor 
product of algebras in $_H{\cal M}$.   
\begin{proposition}
Let $H$ be a quasi-Hopf algebra, $A$ a left $H$-module algebra, 
$B$ an associative algebra, $v:H\rightarrow B$ an algebra map and 
$u:A\rightarrow B^v$ a morphism of left $H$-module algebras. Then we have 
\begin{eqnarray}
&&\mu _{B^v}\circ (u\otimes v)\circ R=\mu _{B^v}\circ (v\otimes u), 
\label{compartic}
\end{eqnarray}
where $R$ is given by (\ref{Rsmash}), that is, the input data for the 
universal property of $A\# H$ is also an input data for the universal 
property of the twisted tensor product $A\otimes _RH_0$ in the 
category $_H{\cal M}$. Moreover, if we denote by $u\# v:A\# H\rightarrow B$  
the algebra map provided by the universal property of $A\# H$ and by 
$w:A\diamond H_0\rightarrow B^v$ the morphism of left $H$-module 
algebras provided by the universal property of $A\diamond H_0=
A\otimes _RH_0$, then we have $(u\# v)\circ \Pi =w$, where $\Pi $ is the 
isomorphism from Proposition \ref{izomleft}. 
\end{proposition}   
\begin{proof}
We denote by $f^{-1}=g^1\otimes g^2=G^1\otimes G^2$ two copies of the element  
$f^{-1}$ given by (\ref{g}). We first check (\ref{compartic}):\\[2mm]
${\;\;\;\;\;\;\;\;\;}$$(\mu _{B^v}\circ (u\otimes v)\circ R)(h\otimes a)$
\begin{eqnarray*}
&=&v(T^1)u(Z^1X^1_1h_1g^1S(Y^2X^2_2)\alpha Y^3X^3\cdot a)
v(S(x^1T^2)\alpha x^2T^3_1)\\
&&v(Z^2X^1_2h_2g^2S(Z^3Y^1X^2_1))v(S(x^3T^3_2))\\
&=&v(T^1Z^1_1X^1_{(1, 1)}h_{(1, 1)}g^1_1S(Y^2X^2_2)_1\alpha _1Y^3_1X^3_1)
u(a)\\
&&v(S(x^1T^2Z^1_2X^1_{(1, 2)}h_{(1, 2)}g^1_2S(Y^2X^2_2)_2\alpha _2Y^3_2X^3_2)\\
&&\alpha x^2T^3_1Z^2X^1_2h_2g^2S(x^3T^3_2Z^3Y^1X^2_1))\\
&{{\rm (\ref{moref}, \;\ref{ca})}\atop =}&
v(T^1Z^1_1X^1_{(1, 1)}h_{(1, 1)}g^1_1G^1
S(W^2t^1_2Y^2_2X^2_{(2, 2)})\alpha W^3t^2Y^3_1X^3_1)u(a)\\
&&v(S(x^1T^2Z^1_2X^1_{(1, 2)}h_{(1, 2)}g^1_2G^2
S(W^1t^1_1Y^2_1X^2_{(2, 1)})\alpha t^3Y^3_2X^3_2)\\
&&\alpha x^2T^3_1Z^2X^1_2h_2g^2S(x^3T^3_2Z^3Y^1X^2_1))\\
&{{\rm (\ref{q3}, \;\ref{q1})}\atop =}&
v(T^1X^1_1h_1Z^1g^1_1G^1
S(W^2t^1_2Y^2_2X^2_{(2, 2)})\alpha W^3t^2Y^3_1X^3_1)u(a)\\
&&v(S(T^2_1X^1_{(2, 1)}h_{(2, 1)}Z^2g^1_2G^2
S(W^1t^1_1Y^2_1X^2_{(2, 1)})\alpha t^3Y^3_2X^3_2)\\
&&\alpha T^2_2X^1_{(2, 2)}h_{(2, 2)}Z^3g^2S(T^3Y^1X^2_1))\\
&{{\rm (\ref{muchmoref})}\atop =}&
v(T^1X^1_1h_1G^1
S(W^2t^1_2Y^2_2X^2_{(2, 2)}Z^3)\alpha W^3t^2Y^3_1X^3_1)u(a)\\
&&v(S(T^2_1X^1_{(2, 1)}h_{(2, 1)}G^2_1g^1
S(W^1t^1_1Y^2_1X^2_{(2, 1)}Z^2)\alpha t^3Y^3_2X^3_2)\\
&&\alpha T^2_2X^1_{(2, 2)}h_{(2, 2)}G^2_2g^2S(T^3Y^1X^2_1Z^1))\\
&{{\rm (\ref{q5})}\atop =}&
v(X^1hS(W^2t^1_2Y^2_2X^2_{(2, 2)}Z^3)\alpha W^3t^2Y^3_1X^3_1)u(a)\\
&&v(S(g^1S(W^1t^1_1Y^2_1X^2_{(2, 1)}Z^2)\alpha t^3Y^3_2X^3_2)
\alpha g^2S(Y^1X^2_1Z^1))\\
&{{\rm (\ref{q3}, \;\ref{q5})}\atop =}&
v(X^1hS(t^1T^2Y^2_2X^2_{(2, 2)}Z^3)\alpha t^2T^3_1Y^3_1X^3_1)u(a)\\
&&v(S(g^1S(T^1Y^2_1X^2_{(2, 1)}Z^2)\alpha t^3T^3_2Y^3_2X^3_2)
\alpha g^2S(Y^1X^2_1Z^1))\\
&{{\rm (\ref{q3}, \;\ref{q1})}\atop =}&
v(X^1hS(t^1T^3_1Y^2X^2_2)\alpha t^2T^3_{(2, 1)}Y^3_1X^3_1)u(a)\\
&&v(S(g^1S(T^2Y^1_2X^2_{(1, 2)})\alpha t^3T^3_{(2, 2)}Y^3_2X^3_2)
\alpha g^2S(T^1Y^1_1X^2_{(1, 1)}))\\
&{{\rm (\ref{q1}, \;\ref{q5})}\atop =}&
v(X^1hS(t^1Y^2X^2_2)\alpha t^2Y^3_1X^3_1)u(a)\\
&&v(S(g^1S(T^2Y^1_2X^2_{(1, 2)})\alpha T^3t^3Y^3_2X^3_2) 
\alpha g^2S(T^1Y^1_1X^2_{(1, 1)}))\\
&{{\rm (\ref{ca})}\atop =}&
v(X^1hS(t^1Y^2X^2_2)\alpha t^2Y^3_1X^3_1)u(a)\\
&&v(S(S(Y^1X^2_1)_1g^1S(T^2)\alpha T^3t^3Y^3_2X^3_2)
\alpha S(Y^1X^2_1)_2g^2S(T^1))\\
&{{\rm (\ref{q5})}\atop =}&
v(X^1hS(t^1X^2)\alpha t^2X^3_1)u(a)
v(S(g^1S(T^2)\alpha T^3t^3X^3_2)
\alpha g^2S(T^1))\\
&{{\rm (\ref{relg}, \;\ref{q6})}\atop =}&
v(X^1hS(t^1X^2)\alpha t^2X^3_1)u(a)v(S(t^3X^3_2))\\
&=&(\mu _{B^v}\circ (v\otimes u))(h\otimes a).
\end{eqnarray*}
Now we check the relation $(u\# v)\circ \Pi =
\mu _{B^v}\circ (u\otimes v)$. We compute:
\begin{eqnarray*}
((u\# v)\circ \Pi )(a\otimes h)&=&(u\# v)(x^1\cdot a\# x^2hS(x^3))\\
&{{\rm (\ref{unique})}\atop =}&v(X^1)u(x^1\cdot a)v(S(X^2)\alpha X^3
x^2hS(x^3))\\
&=&v(X^1x^1_1)u(a)v(S(X^2x^1_2)\alpha X^3x^2hS(x^3))\\
&{{\rm (\ref{q3}, \;\ref{q5})}\atop =}&
v(X^1)u(a)v(S(x^1X^2)\alpha x^2X^3_1hS(x^3X^3_2))\\
&=&(\mu _{B^v}\circ (u\otimes v))(a\otimes h),
\end{eqnarray*}
finishing the proof.
\end{proof}
${\;\;\;}$Let $H$ be an ordinary Hopf algebra and $A$ a left $H$-module 
algebra. We can regard $A\# H$ as a twisted tensor product between the 
associative algebras $A$ and $H$, and thus it has a corresponding universal 
property, which is easily seen to be equivalent to the usual universal 
property of $A\# H$. \\
${\;\;\;}$We would like to have a quasi-Hopf analogue of this situation, 
that is, a universal property of the quasi-Hopf smash product $A\# H$ 
resembling the universal property of a twisted tensor product. It is not clear 
{\it a priori} how to formulate such a universal property, since $A\# H$ 
itself is not a twisted tensor product, neither of associative algebras nor 
of quasialgebras.\\
${\;\;\;}$Our strategy to arrive at such a formulation is the following. 
We obtain first a universal property of the generalized diagonal crossed 
product ${\cal A}\bowtie {\mb A}$, which generalizes (and is inspired by) 
both the universal property of the diagonal crossed product 
$H^*\bowtie {\mb A}$ over a finite dimensional quasi-Hopf algebra $H$ from 
\cite{hn1} and the universal property of a generalized diagonal crossed 
product over a Hopf algebra, regarded as a twisted tensor product 
of associative algebras. Then, since $A\# H$ is a particular case of 
${\cal A}\bowtie {\mb A}$, we will obtain the desired universal property 
of $A\# H$.\\
${\;\;\;}$If $H$ is a quasi-Hopf algebra, ${\cal A}$ an $H$-bimodule algebra 
and ${\mb A}$ an $H$-bicomodule algebra, we denote by 
$\Gamma :{\cal A}\rightarrow {\cal A}\bowtie {\mb A}$ and 
$j:{\mb A}\rightarrow {\cal A}\bowtie {\mb A}$ the linear maps given by 
\begin{eqnarray*}
&&\Gamma (\varphi )=(\tilde{p}^1_{\rho })_{[-1]}\cdot \varphi \cdot 
S^{-1}(\tilde{p}^2_{\rho })\bowtie (\tilde{p}^1_{\rho })_{[0]},\;\;\;
j(u)=1_{{\cal A}}\bowtie u,
\end{eqnarray*}
for all $\varphi \in {\cal A}$, $u\in {\mb A}$, where 
$\tilde{p}^1_{\rho }\otimes \tilde{p}^2_{\rho }$ is given by (\ref{tpqr}). 
By \cite{bpvo}, the map $\Gamma $ has the property that 
${\cal A}\bowtie {\mb A}$ is generated as algebra by ${\mb A}$ and 
$\Gamma ({\cal A})$, while $j$ is obviously an algebra map.  
\begin{proposition} \label{univdiagonal}
Let $H$ be a quasi-Hopf algebra, ${\cal A}$ an $H$-bimodule algebra, 
${\mb A}$ an $H$-bicomodule algebra, $B$ an associative algebra, 
$\gamma :{\mb A}\rightarrow B$ an algebra map and $v:{\cal A}\rightarrow B$ 
a linear map such that the following conditions are satisfied:
\begin{eqnarray}
&&\gamma (u_{<0>})v(\varphi \cdot u_{<1>})=v(u_{[-1]}\cdot \varphi )
\gamma (u_{[0]}), \label{cond1} \\
&&v(\varphi \varphi ')=\gamma (\tilde{X}^1_{\rho })v(\theta ^1
\tilde{X}^1_{\lambda }\cdot \varphi \cdot \tilde{X}^2_{\rho })
\gamma (\theta ^2)v(\tilde{X}^2_{\lambda }\cdot \varphi '\cdot 
\tilde{X}^3_{\rho }\theta ^3)\gamma (\tilde{X}^3_{\lambda }), 
\label{cond2} \\
&&v(1_{{\cal A}})=1_B, \label{cond3}
\end{eqnarray}  
for all $\varphi , \varphi '\in {\cal A}$ and $u\in {\mb A}$. Then there 
exists a unique algebra map $w:{\cal A}\bowtie {\mb A}\rightarrow B$ 
such that $w\circ \Gamma =v$ and $w\circ j=\gamma $. Moreover, 
$w$ is given by the formula 
\begin{eqnarray}
&&w(\varphi \bowtie u)=\gamma (\tilde{q}^1_{\rho })v(\varphi \cdot 
\tilde{q}^2_{\rho })\gamma (u), \label{unica}
\end{eqnarray}
for all $\varphi \in {\cal A}$ and $u\in {\mb A}$, where 
$\tilde{q}_{\rho }=\tilde{q}^1_{\rho }\otimes \tilde{q}^2_{\rho }$ is the 
element given by the formula (\ref{tpqr}). 
\end{proposition}
\begin{proof}
We first prove the uniqueness of $w$. By \cite{bpvo}, the map $\Gamma $ 
satisfies the relation 
\begin{eqnarray*}
&&\varphi \bowtie u=(1_{{\cal A}}\bowtie \tilde{q}^1_{\rho })
\Gamma (\varphi \cdot \tilde{q}^2_{\rho })(1_{{\cal A}}\bowtie u),
\end{eqnarray*}
for all $\varphi \in {\cal A}$ and $u\in {\mb A}$, hence we can write
\begin{eqnarray*}
w(\varphi \bowtie u)&=&w(j(\tilde{q}^1_{\rho })\Gamma (\varphi \cdot 
\tilde{q}^2_{\rho })j(u))\\
&=&w(j(\tilde{q}^1_{\rho }))w(\Gamma (\varphi \cdot 
\tilde{q}^2_{\rho }))w(j(u))\\
&=&\gamma (\tilde{q}^1_{\rho })v(\varphi \cdot 
\tilde{q}^2_{\rho })\gamma (u),
\end{eqnarray*}
showing that $w$ is indeed unique. We now prove the existence part. 
Define $w$ by formula (\ref{unica}); it is obvious that $w$ is unital and 
satisfies $w\circ j=\gamma $. We check now that $w\circ \Gamma =v$:
\begin{eqnarray*}
(w\circ \Gamma )(\varphi )&=&w((\tilde{p}^1_{\rho })_{[-1]}\cdot \varphi 
\cdot S^{-1}(\tilde{p}^2_{\rho })\bowtie (\tilde{p}^1_{\rho })_{[0]})\\ 
&=&\gamma (\tilde{q}^1_{\rho })v((\tilde{p}^1_{\rho })_{[-1]}\cdot \varphi 
\cdot S^{-1}(\tilde{p}^2_{\rho })\tilde{q}^2_{\rho })
\gamma ((\tilde{p}^1_{\rho })_{[0]})\\
&{{\rm (\ref{cond1})}\atop =}&
\gamma (\tilde{q}^1_{\rho })\gamma ((\tilde{p}^1_{\rho })_{<0>})v(\varphi  
\cdot S^{-1}(\tilde{p}^2_{\rho })\tilde{q}^2_{\rho }
(\tilde{p}^1_{\rho })_{<1>})\\
&{{\rm (\ref{tpqr2a})}\atop =}&v(\varphi ).
\end{eqnarray*} 
Thus, the only thing left to prove is that $w$ is multiplicative. We denote 
by $\tilde{Q}^1_{\rho }\otimes \tilde{Q}^2_{\rho }$ another copy of the 
element $\tilde{q}_{\rho }$, and we record the obvious relation 
\begin{eqnarray}
&&\tilde{Q}^1_{\rho }\tilde{x}^1_{\rho }\otimes S^{-1}(\tilde{x}^3_{\rho })
\tilde{Q}^2_{\rho }\tilde{x}^2_{\rho }=1\otimes S^{-1}(\alpha ). 
\label{simplu}
\end{eqnarray}  
Now we compute:\\[2mm]
${\;\;}$$w((\varphi \bowtie u)(\varphi '\bowtie u'))$
\begin{eqnarray*}
&=&\gamma (\tilde{q}^1_{\rho })v([\Omega ^1\cdot \varphi \cdot \Omega ^5
(\tilde{q}^2_{\rho })_1][\Omega ^2u_{<0>_{[-1]}}\cdot \varphi '\cdot 
S^{-1}(u_{<1>})\Omega ^4(\tilde{q}^2_{\rho })_2])\gamma (\Omega ^3
u_{<0>_{[0]}}u')\\
&{{\rm (\ref{cond2})}\atop =}&
\gamma (\tilde{q}^1_{\rho })\gamma (\tilde{X}^1_{\rho })
v(\theta ^1\tilde{X}^1_{\lambda }\Omega ^1\cdot \varphi \cdot \Omega ^5
(\tilde{q}^2_{\rho })_1\tilde{X}^2_{\rho })\gamma (\theta ^2)\\
&&v(\tilde{X}^2_{\lambda }
\Omega ^2u_{<0>_{[-1]}}\cdot \varphi '\cdot 
S^{-1}(u_{<1>})\Omega ^4(\tilde{q}^2_{\rho })_2\tilde{X}^3_{\rho }
\theta ^3)\gamma (\tilde{X}^3_{\lambda })
\gamma (\Omega ^3u_{<0>_{[0]}}u')\\
&{{\rm (\ref{formomega})}\atop =}&
\gamma (\tilde{q}^1_{\rho }\tilde{X}^1_{\rho })
v(\theta ^1(\tilde{Y}^1_{\rho })_{[-1]}\overline{\theta }^1\cdot \varphi 
\cdot S^{-1}(f^2\tilde{Y}^3_{\rho })
(\tilde{q}^2_{\rho })_1\tilde{X}^2_{\rho })\gamma (\theta ^2)\\
&&v((\tilde{Y}^1_{\rho })_{[0]_{[-1]}}\overline{\theta }^2_{[-1]}
u_{<0>_{[-1]}}\cdot \varphi '\cdot  
S^{-1}(f^1\tilde{Y}^2_{\rho }\overline{\theta }^3u_{<1>})
(\tilde{q}^2_{\rho })_2\tilde{X}^3_{\rho }\theta ^3)\\
&&\gamma ((\tilde{Y}^1_{\rho })_{[0]_{[0]}}\overline{\theta }^2_{[0]}
u_{<0>_{[0]}})\gamma (u')\\
&{{\rm (\ref{cond1})}\atop =}&
\gamma (\tilde{q}^1_{\rho }\tilde{X}^1_{\rho })
v(\theta ^1(\tilde{Y}^1_{\rho })_{[-1]}\overline{\theta }^1\cdot \varphi 
\cdot S^{-1}(f^2\tilde{Y}^3_{\rho })
(\tilde{q}^2_{\rho })_1\tilde{X}^2_{\rho })
\gamma (\theta ^2(\tilde{Y}^1_{\rho })_{[0]_{<0>}}\overline{\theta }^2_{<0>}
u_{<0>_{<0>}})\\
&&v(\varphi '\cdot  
S^{-1}(f^1\tilde{Y}^2_{\rho }\overline{\theta }^3u_{<1>})
(\tilde{q}^2_{\rho })_2\tilde{X}^3_{\rho }\theta ^3
(\tilde{Y}^1_{\rho })_{[0]_{<1>}}\overline{\theta }^2_{<1>}
u_{<0>_{<1>}})\gamma (u')\\
&{{\rm (\ref{bca1})}\atop =}&
\gamma (\tilde{q}^1_{\rho }\tilde{X}^1_{\rho })
v((\tilde{Y}^1_{\rho })_{<0>_{[-1]}}\theta ^1\overline{\theta }^1
\cdot \varphi \cdot S^{-1}(f^2\tilde{Y}^3_{\rho })
(\tilde{q}^2_{\rho })_1\tilde{X}^2_{\rho })
\gamma ((\tilde{Y}^1_{\rho })_{<0>_{[0]}}\theta ^2\overline{\theta }^2_{<0>} 
u_{<0>_{<0>}})\\
&&v(\varphi '\cdot  
S^{-1}(f^1\tilde{Y}^2_{\rho }\overline{\theta }^3u_{<1>})
(\tilde{q}^2_{\rho })_2\tilde{X}^3_{\rho }
(\tilde{Y}^1_{\rho })_{<1>}\theta ^3\overline{\theta }^2_{<1>}
u_{<0>_{<1>}})\gamma (u')\\
&{{\rm (\ref{bca3})}\atop =}&
\gamma (\tilde{q}^1_{\rho }\tilde{X}^1_{\rho })
v((\tilde{Y}^1_{\rho })_{<0>_{[-1]}}(\tilde{y}^1_{\rho })_{[-1]}\theta ^1
\cdot \varphi \cdot S^{-1}(f^2\tilde{Y}^3_{\rho })
(\tilde{q}^2_{\rho })_1\tilde{X}^2_{\rho })\\
&&\gamma ((\tilde{Y}^1_{\rho })_{<0>_{[0]}}(\tilde{y}^1_{\rho })_{[0]}
\theta ^2\tilde{Z}^1_{\rho }u_{<0>_{<0>}})\\
&&v(\varphi '\cdot  
S^{-1}(f^1\tilde{Y}^2_{\rho }\tilde{y}^3_{\rho }\theta ^3_2
\tilde{Z}^3_{\rho }u_{<1>})
(\tilde{q}^2_{\rho })_2\tilde{X}^3_{\rho }
(\tilde{Y}^1_{\rho })_{<1>}\tilde{y}^2_{\rho }\theta ^3_1\tilde{Z}^2_{\rho }
u_{<0>_{<1>}})\gamma (u')\\
&{{\rm (\ref{cond1})}\atop =}&
\gamma (\tilde{q}^1_{\rho }\tilde{X}^1_{\rho }
(\tilde{Y}^1_{\rho })_{<0>_{<0>}}(\tilde{y}^1_{\rho })_{<0>})v(\theta ^1 
\cdot \varphi \cdot S^{-1}(f^2\tilde{Y}^3_{\rho })
(\tilde{q}^2_{\rho })_1\tilde{X}^2_{\rho }
(\tilde{Y}^1_{\rho })_{<0>_{<1>}}(\tilde{y}^1_{\rho })_{<1>})\\
&&\gamma (\theta ^2\tilde{Z}^1_{\rho }u_{<0>_{<0>}})\\
&&v(\varphi '\cdot   
S^{-1}(f^1\tilde{Y}^2_{\rho }\tilde{y}^3_{\rho }\theta ^3_2
\tilde{Z}^3_{\rho }u_{<1>})
(\tilde{q}^2_{\rho })_2\tilde{X}^3_{\rho }
(\tilde{Y}^1_{\rho })_{<1>}\tilde{y}^2_{\rho }\theta ^3_1\tilde{Z}^2_{\rho }
u_{<0>_{<1>}})\gamma (u')\\
&{{\rm (\ref{rca1}, \;\ref{rca2})}\atop =}&
\gamma (\tilde{q}^1_{\rho }
(\tilde{Y}^1_{\rho })_{<0>}\tilde{y}^1_{\rho }\tilde{X}^1_{\rho })
v(\theta ^1\cdot \varphi \cdot S^{-1}(f^2\tilde{Y}^3_{\rho })
(\tilde{q}^2_{\rho })_1
(\tilde{Y}^1_{\rho })_{<1>_1}(\tilde{y}^2_{\rho })_1x^1\tilde{X}^2_{\rho })\\
&&\gamma (\theta ^2\tilde{Z}^1_{\rho }u_{<0>_{<0>}})
v(\varphi '\cdot   
S^{-1}(f^1\tilde{Y}^2_{\rho }\tilde{y}^3_{\rho }x^3
(\tilde{X}^3_{\rho })_2\theta ^3_2\tilde{Z}^3_{\rho }u_{<1>})\\
&&(\tilde{q}^2_{\rho })_2
(\tilde{Y}^1_{\rho })_{<1>_2}(\tilde{y}^2_{\rho })_2x^2
(\tilde{X}^3_{\rho })_1\theta ^3_1\tilde{Z}^2_{\rho }
u_{<0>_{<1>}})\gamma (u')\\
&{{\rm (\ref{rca1}, \;\ref{relfoartemare})}\atop =}&
\gamma (\tilde{q}^1_{\rho }(\tilde{Q}^1_{\rho })_{<0>}\tilde{x}^1_{\rho }
\tilde{y}^1_{\rho }\tilde{X}^1_{\rho })
v(\theta ^1\cdot \varphi \cdot \tilde{q}^2_{\rho }(\tilde{Q}^1_{\rho })_{<1>} 
\tilde{x}^2_{\rho }(\tilde{y}^2_{\rho })_1x^1\tilde{X}^2_{\rho })
\gamma (\theta ^2u_{<0>}\tilde{Z}^1_{\rho })\\
&&v(\varphi '\cdot S^{-1}(\tilde{y}^3_{\rho }x^3
(\tilde{X}^3_{\rho })_2\theta ^3_2u_{<1>_2}\tilde{Z}^3_{\rho })
\tilde{Q}^2_{\rho }\tilde{x}^3_{\rho }(\tilde{y}^2_{\rho })_2x^2
(\tilde{X}^3_{\rho })_1\theta ^3_1u_{<1>_1}\tilde{Z}^2_{\rho })\gamma (u')\\
&{{\rm (\ref{rca2})}\atop =}&
\gamma (\tilde{q}^1_{\rho }(\tilde{Q}^1_{\rho })_{<0>}
(\tilde{x}^1_{\rho })_{<0>})
v(\theta ^1\cdot \varphi \cdot \tilde{q}^2_{\rho }(\tilde{Q}^1_{\rho })_{<1>} 
(\tilde{x}^1_{\rho })_{<1>}) 
\gamma (\theta ^2u_{<0>}\tilde{Z}^1_{\rho })\\
&&v(\varphi '\cdot S^{-1}(\tilde{x}^3_{\rho }
\theta ^3_2u_{<1>_2}\tilde{Z}^3_{\rho })
\tilde{Q}^2_{\rho }\tilde{x}^2_{\rho }
\theta ^3_1u_{<1>_1}\tilde{Z}^2_{\rho })\gamma (u')\\
&{{\rm (\ref{simplu}, \;\ref{q5})}\atop =}&
\gamma (\tilde{q}^1_{\rho })
v(\varphi \cdot \tilde{q}^2_{\rho }) 
\gamma (u\tilde{Z}^1_{\rho })
v(\varphi '\cdot S^{-1}(\tilde{Z}^3_{\rho })S^{-1}(\alpha )
\tilde{Z}^2_{\rho })\gamma (u')\\
&=&\gamma (\tilde{q}^1_{\rho })
v(\varphi \cdot \tilde{q}^2_{\rho }) 
\gamma (u)\gamma (\tilde{Q}^1_{\rho })
v(\varphi '\cdot \tilde{Q}^2_{\rho })\gamma (u')\\
&=&w(\varphi \bowtie u)w(\varphi '\bowtie u'),
\end{eqnarray*} 
finishing the proof.
\end{proof}
\begin{remark}{\em 
This universal property of ${\cal A}\bowtie {\mb A}$ looks quite 
different from the universal property of $H^*\bowtie {\mb A}$ given 
in \cite{hn1}, Theorem II, which is formulated in terms of so-called 
{\it normal coherent intertwiners}, but for finite 
dimensional $H$ and ${\cal A}=H^*$ the two universal properties are actually 
equivalent. Indeed, if $T=T^1\otimes T^2\in H\otimes B$ is a 
normal coherent intertwiner as in \cite{hn1}, then we can define the map 
$v:H^*\rightarrow B$, $v(\varphi )=\varphi (T^1)T^2$, and the three  
conditions (\ref{cond1}), (\ref{cond2}) and (\ref{cond3}) for $v$ 
follow from the three conditions in \cite{hn1} defining  
a normal coherent intertwiner.}   
\end{remark}
${\;\;\;}$We intend now to prove that the algebra isomorphism 
$\nu :{\cal A}\bowtie {\mb A}\rightarrow {\cal A}\nat {\mb A}$ given by 
(\ref{formulanu}) may be naturally reobtained by using 
Proposition \ref{univdiagonal}. We define the map 
\begin{eqnarray}
&&\Lambda :{\cal A}\rightarrow {\cal A}\nat {\mb A}, \;\;\;
\Lambda (\varphi )=\theta ^1\cdot \varphi \cdot \theta ^3\nat \theta ^2.
\label{formulaLambda} 
\end{eqnarray}
We claim that the conditions in Proposition \ref{univdiagonal} 
are satisfied for $B={\cal A}\nat {\mb A}$, $\gamma =j$ and $v=\Lambda $. 
Indeed, condition (\ref{cond1}) follows easily from (\ref{bca1}) and 
condition (\ref{cond3}) is trivial, so we only have to check (\ref{cond2}). 
We denote by $\tilde{\theta }^1\otimes \tilde{\theta }^2\otimes 
\tilde{\theta }^3$ one more copy of $\Phi _{\lambda , \rho }^{-1}$ and we  
compute:\\[2mm]
${\;\;}$$j(\tilde{X}^1_{\rho })\Lambda (\theta ^1\tilde{X}^1_{\lambda }
\cdot \varphi \cdot \tilde{X}^2_{\rho })j(\theta ^2)
\Lambda (\tilde{X}^2_{\lambda }\cdot \varphi '\cdot \tilde{X}^3_{\rho }
\theta ^3)j(\tilde{X}^3_{\lambda })$
\begin{eqnarray*}
&=&((\tilde{X}^1_{\rho })_{[-1]}\overline{\theta }^1\theta ^1
\tilde{X}^1_{\lambda }\cdot \varphi \cdot \tilde{X}^2_{\rho }
\overline{\theta }^3\theta ^2_{<1>}\nat (\tilde{X}^1_{\rho })_{[0]}
\overline{\theta }^2\theta ^2_{<0>})\\
&&(\tilde{\theta }^1\tilde{X}^2_{\lambda }\cdot \varphi '\cdot 
\tilde{X}^3_{\rho }\theta ^3\tilde{\theta }^3(\tilde{X}^3_{\lambda })_{<1>}
\nat \tilde{\theta }^2(\tilde{X}^3_{\lambda })_{<0>})\\
&{{\rm (\ref{bca3})}\atop =}&
(\theta ^1\tilde{X}^1_{\lambda }\cdot \varphi \cdot \theta ^3_1
\tilde{Y}^2_{\rho }\nat \theta ^2\tilde{Y}^1_{\rho })
(\tilde{\theta }^1\tilde{X}^2_{\lambda }\cdot \varphi '\cdot 
\theta ^3_2\tilde{Y}^3_{\rho }\tilde{\theta }^3(\tilde{X}^3_{\lambda })_{<1>}
\nat \tilde{\theta }^2(\tilde{X}^3_{\lambda })_{<0>})\\
&=&(\tilde{x}^1_{\lambda }\theta ^1\tilde{X}^1_{\lambda }\cdot \varphi 
\cdot \theta ^3_1\tilde{Y}^2_{\rho }\overline{\theta }^3
\tilde{\theta }^2_{<1>}(\tilde{X}^3_{\lambda })_{<0>_{<1>}}
\tilde{x}^2_{\rho })\\
&&(\tilde{x}^2_{\lambda }\theta ^2_{[-1]}
(\tilde{Y}^1_{\rho })_{[-1]}\overline{\theta }^1\tilde{\theta }^1
\tilde{X}^2_{\lambda }\cdot \varphi '\cdot 
\theta ^3_2\tilde{Y}^3_{\rho }\tilde{\theta }^3
(\tilde{X}^3_{\lambda })_{<1>}\tilde{x}^3_{\rho })
\nat \tilde{x}^3_{\lambda }\theta ^2_{[0]}(\tilde{Y}^1_{\rho })_{[0]}
\overline{\theta }^2\tilde{\theta }^2_{<0>}
(\tilde{X}^3_{\lambda })_{<0>_{<0>}}\tilde{x}^1_{\rho }\\
&{{\rm (\ref{bca3})}\atop =}&
(\tilde{x}^1_{\lambda }\theta ^1\tilde{X}^1_{\lambda }\cdot \varphi 
\cdot \theta ^3_1\overline{\theta }^3_1\tilde{Y}^2_{\rho }
(\tilde{X}^3_{\lambda })_{<0>_{<1>}}
\tilde{x}^2_{\rho })\\
&&(\tilde{x}^2_{\lambda }\theta ^2_{[-1]}\overline{\theta }^1
\tilde{X}^2_{\lambda }\cdot \varphi '\cdot 
\theta ^3_2\overline{\theta }^3_2\tilde{Y}^3_{\rho }
(\tilde{X}^3_{\lambda })_{<1>}\tilde{x}^3_{\rho })
\nat \tilde{x}^3_{\lambda }\theta ^2_{[0]}
\overline{\theta }^2\tilde{Y}^1_{\rho }
(\tilde{X}^3_{\lambda })_{<0>_{<0>}}\tilde{x}^1_{\rho }\\
&{{\rm (\ref{bca2})}\atop =}&
(\theta ^1_1\tilde{x}^1_{\lambda }\tilde{X}^1_{\lambda }\cdot \varphi 
\cdot \theta ^3_1(\tilde{x}^3_{\lambda })_{<1>_1}\tilde{Y}^2_{\rho }
(\tilde{X}^3_{\lambda })_{<0>_{<1>}}
\tilde{x}^2_{\rho })\\
&&(\theta ^1_2\tilde{x}^2_{\lambda }
\tilde{X}^2_{\lambda }\cdot \varphi '\cdot 
\theta ^3_2(\tilde{x}^3_{\lambda })_{<1>_2}\tilde{Y}^3_{\rho }
(\tilde{X}^3_{\lambda })_{<1>}\tilde{x}^3_{\rho })
\nat \theta ^2(\tilde{x}^3_{\lambda })_{<0>}\tilde{Y}^1_{\rho }
(\tilde{X}^3_{\lambda })_{<0>_{<0>}}\tilde{x}^1_{\rho }\\
&{{\rm (\ref{rca1})}\atop =}&
(\theta ^1_1\cdot \varphi \cdot \theta ^3_1)
(\theta ^1_2\cdot \varphi '\cdot \theta ^3_2)
\nat \theta ^2\\
&=&\Lambda (\varphi \varphi '),
\end{eqnarray*}
finishing the proof of (\ref{cond2}). We can thus apply 
Proposition \ref{univdiagonal}, obtaining  
an algebra map $w:{\cal A}\bowtie {\mb A}\rightarrow {\cal A}\nat {\mb A}$ 
satisfying $w\circ \Gamma =\Lambda $ and $w\circ j=j$, and given by the 
formula 
\begin{eqnarray*}
w(\varphi \bowtie u)&=&j(\tilde{q}^1_{\rho })\Lambda (\varphi \cdot 
\tilde{q}^2_{\rho })j(u)\\
&=&(\tilde{q}^1_{\rho })_{[-1]}\theta ^1\cdot \varphi \cdot  
\tilde{q}^2_{\rho }\theta ^3u_{<1>}\nat (\tilde{q}^1_{\rho })_{[0]}
\theta ^2u_{<0>}.
\end{eqnarray*}
On the other hand, by using (\ref{bca3}) and (\ref{q5}), one can easily 
check the following relation: 
\begin{eqnarray*}
&&(\tilde{q}^1_{\rho })_{[-1]}\theta ^1\otimes 
(\tilde{q}^1_{\rho })_{[0]}\theta ^2\otimes \tilde{q}^2_{\rho }\theta ^3=
\Theta ^1\otimes \tilde{q}^1_{\rho }\Theta ^2_{<0>}\otimes 
S^{-1}(\Theta ^3)\tilde{q}^2_{\rho }\Theta ^2_{<1>},
\end{eqnarray*}
which implies that the map $w$ above coincides with the algebra 
isomorphism $\nu $ given by (\ref{formulanu}).\\[2mm]      
${\;\;\;}$As a consequence of Proposition \ref{univdiagonal}, we  
immediately obtain the following new kind of universal property for the 
quasi-Hopf smash product:
\begin{proposition} \label{newunivsmash}
Let $H$ be a quasi-Hopf algebra and $A$ a left $H$-module algebra. 
Denote by $i:A\rightarrow A\# H$, $i(a)=a\# 1$ and $j:H\rightarrow A\# H$, 
$j(h)=1\# h$. Let $B$ be an associative algebra,  
$\gamma :H\rightarrow B$ an algebra map and $v:A\rightarrow B$  
a linear map satisfying the following conditions:
\begin{eqnarray}
&&\gamma (h)v(a)=v(h_1\cdot a)
\gamma (h_2), \label{partcond1} \\
&&v(aa')=v(X^1\cdot a)v(X^2\cdot a')\gamma (X^3), 
\label{partcond2} \\
&&v(1_{A})=1_B, \label{partcond3}
\end{eqnarray}  
for $a, a'\in A$, $h\in H$. Then there  
exists a unique algebra map $w:A\# H\rightarrow B$  
such that $w\circ i=v$ and $w\circ j=\gamma $. Moreover, 
$w$ is given by the formula 
$w(a\# h)=v(a)\gamma (h)$, 
for $a\in A$, $h\in H$.  
\end{proposition}
${\;\;\;}$The relation between the two universal properties of $A\# H$ 
is the following: Proposition \ref{newunivsmash} implies the usual universal 
property (unlike the case of ordinary Hopf algebras, the converse does 
not seem to hold). Indeed, if $B$ is an associative algebra, 
$\gamma :H\rightarrow B$ an algebra map and $u:A\rightarrow B^{\gamma }$ a 
morphism of left $H$-module algebras, define the map 
\begin{eqnarray*}
&&v:A\rightarrow B, \;\;\;v(a)=\gamma (q^1)u(a)\gamma (S(q^2)), 
\end{eqnarray*}
where $q_R=q^1\otimes q^2$ is given by (\ref{psiq}). Then one can check that  
the maps $\gamma $ and $v$ satisfy the hypotheses of Proposition 
\ref{newunivsmash}, and obviously the unique algebra map 
$w:A\# H\rightarrow B$ provided by Proposition \ref{newunivsmash} 
coincides with the map $u\# \gamma $ given by (\ref{unique}).\\
${\;\;\;}$Proposition \ref{newunivsmash} may be easily extended to a 
universal property of the two-sided smash product:
\begin{proposition} 
Let $H$ be a quasi-Hopf algebra, $A$ a left $H$-module algebra and $B$ 
a right $H$-module algebra.  
Denote by $i_A$, $i_B$, $j$ the standard inclusions of $A$, $B$, $H$ 
respectively into $A\# H\# B$. Let $X$ be an associative algebra,   
$\gamma :H\rightarrow X$ an algebra map and $v_A:A\rightarrow X$, 
$v_B:B\rightarrow X$ two linear maps satisfying the following conditions:
\begin{eqnarray*}
&&\gamma (h)v_A(a)=v_A(h_1\cdot a)
\gamma (h_2), \;\;\;
v_A(aa')=v_A(X^1\cdot a)v_A(X^2\cdot a')\gamma (X^3), \\
&&v_B(b)\gamma (h)=\gamma (h_1)v_B(b\cdot h_2), \;\;\;
v_B(bb')=\gamma (X^1)v_B(b\cdot X^2)v_B(b'\cdot X^3), \\
&&v_A(1_A)=1_X=v_B(1_B), \;\;\;
v_B(b)v_A(a)=v_A(x^1\cdot a)\gamma (x^2)v_B(b\cdot x^3), 
\end{eqnarray*}  
for all $a, a'\in A$, $b, b'\in B$, $h\in H$. Then there   
exists a unique algebra map $w:A\# H\# B\rightarrow X$   
such that $w\circ i_A=v_A$, $w\circ i_B=v_B$ and $w\circ j=\gamma $. 
Moreover, $w$ is given by the formula $w(a\# h\# b)=v_A(a)\gamma (h)v_B(b)$,  
for all $a\in A$, $h\in H$, $b\in B$.  
\end{proposition}
${\;\;\;}$For instance, the algebra isomorphism $A\# H\# B\simeq 
(A\otimes B)\nat H$ from \cite{pv} may be easily reobtained by using this 
universal property.  

\end{document}